\documentclass[fleqn,10pt]{wlscirep}
\usepackage[utf8]{inputenc}
\usepackage[T1]{fontenc}
\usepackage[ruled,vlined] {algorithm2e}
\usepackage{mathtools}
\usepackage[table]{xcolor}
\usepackage{array}

\usepackage{extpfeil}
\usepackage{float}

\newcommand{\candLibMatrix}{\Theta}

\title{SODAs: Sparse Optimization for the Discovery of Differential and Algebraic Equations}

\author[1]{Manu Jayadharan}
\author[1]{Christina Catlett}
\author[2]{Arthur N. Montanari}
\author[1,*]{Niall M. Mangan}
\affil[1]{Department of Engineering Sciences and Applied Mathematics, Northwestern University, Evanston, IL 60208, USA}
\affil[2]{Department of Physics and Astronomy, Northwestern University, Evanston, IL 60208, USA}

\affil[*]{niall.mangan@northwestern.edu}

\begin{abstract}
Differential-algebraic equations (DAEs) integrate ordinary differential equations (ODEs) with algebraic constraints, providing a fundamental framework for developing models of dynamical systems characterized by timescale separation, conservation laws, and physical constraints. While sparse optimization has revolutionized model development by allowing data-driven discovery of parsimonious models from a library of possible equations, existing approaches for dynamical systems assume DAEs can be reduced to ODEs by eliminating variables before model discovery. This assumption limits the applicability of such methods for DAE systems with unknown constraints and time scales. 
We introduce Sparse Optimization for Differential-Algebraic Systems (SODAs), a data-driven method for the identification of DAEs in their explicit form. By discovering the algebraic and dynamic components sequentially without prior identification of the algebraic variables, this approach leads to a sequence of convex optimization problems. It has the advantage of discovering interpretable models that preserve the structure of the underlying physical system. To this end, SODAs improves numerical stability when handling high correlations between library terms, caused by near-perfect algebraic relationships, by iteratively refining the conditioning of the candidate library.
We demonstrate the performance of our method on biological, mechanical, and electrical systems, showcasing its robustness to noise in both simulated time series and real-time experimental data. 
\end{abstract}
\begin{document}

\flushbottom
\maketitle
%
%

\section{Introduction}\label{Sec:Intro}
Data-driven discovery of dynamical models of complex systems promises to rapidly accelerate our understanding, from learning new equations governing cellular biomechanics [\citen{SCHMITT2024481}] to the automatic generation of control models for robotics [\citen{RoboticsAmini,RoboticsTaylor}].  A major challenge of many complex dynamical systems is the separation of timescales and the integration of physical constraints and conservation laws. Timescale separation is ubiquitous in complex biological and physical systems, from metabolic networks to multi-body robotics and power grids. Even in well-studied systems such as those following mass-action kinetics and deterministic Newtonian physics, numerical simulation and parameter estimation are challenging due to the inherent stiffness of such systems [\citen{stiff_ODE_1,Stiff_ODE_2,byrne1987stiff}]. In this context, model selection is an even greater challenge as the structural form of the model, including its conservation laws and relevant timescales, may be unknown. In this work, we propose a method for data-driven discovery of complex systems by designing sparse-model selection within a differential-algebraic equation (DAE) framework, where conservation laws and quasi-steady-state equations are learned from time-series data as algebraic equations.

Model identification or selection focuses on the general challenge of discovering the functional form of a model using data. For many dynamical systems, including those following mass-action kinetics or electromechanical conservation laws, it is reasonable to enumerate a library of possible interaction terms or candidate functions \(\candLibMatrix = \{ \theta_1(\mathbf{x}), \ldots, \theta_J(\mathbf{x}) \}\), each representing a meaningful physical quantity as a function of the system state $\mathbf{x}(t)$. The key challenge is to downselect the most suitable subset of terms from this library. For example, a system of differential equations can often be described using a linear combination of nonlinear functions: $\dot{x}_i = \sum_{j=1}^{J} p_j^i \theta_j(\mathbf{x}),$ for each state variable $i = 1, \ldots, d$, where $p_i^j$ are parameters. Popular candidates for \(\candLibMatrix\) include monomial bases, rational functions, and Fourier bases. There exist several methods to identify nonzero coefficients $p_j^i \neq 0$, thereby learning which functions $\theta_j(\mathbf{x})$ influence the dynamics of state $x_i$, including Bayesian variable selection [\citen{Bayesian_Mitchell}], best subset selection [\citen{subset_sel_Miller,best_subset}], and sparse regression [\citen{sindy}]. Of these, sparse regression has gained popularity due to its relative simplicity and low computational cost. Among the best known examples, Sparse Optimization for Nonlinear Dynamics (SINDy) finds a parsimonious model by minimizing 
\begin{align}
    L_i(\mathbf{p}^i) &= \sum_{r=1}^{N} \| \dot{x}_i(t_r) - \sum_{j=1}^{J} p_j^i \theta_j\left(x_1(t_r), x_2(t_r), \ldots, x_d(t_r) \right) \|^2 + \alpha \|\mathbf{p}^i\|_0, \quad i = 1, \ldots, d, \label{eq:loss-cont}
\end{align}
where \(\alpha\) controls the weight of the penalty term enforcing sparsity and $N$ is the number of data samples, and $\mathbf{p}^i$ is a vector of parameters for the $i^{th}$ state equation. If all states $x_i$ are measured and the derivatives $\dot{x}_i$ can be measured or numerically estimated, then the model is linear with respect to the unknown parameters and hence minimization can be performed through alternating least-squares and hard-thresholding or regularized approaches such as LASSO [\citen{Lasso_tibshirani}], Elastic-net  [\citen{EN_orig}], or SR3 [\citen{SR3_orig}]. Some of these methods replace the $l_0$ constraint in \eqref{eq:loss-cont} with either an $l_1$ or $l_2$ penalty, which leads to a convex optimization problem and is often combined with sequential thresholding. Each of these methods has its own advantages and disadvantages. Previous work, [\citen{Pathak2024OnTD}],  has shown the sub-optimality of LASSO compared to the soft-thresholded ordinary least squares estimator (STOLS). Moreover, [\citen{zheng2014high}] shows that hard thresholding combined with $l_2$ shrinkage can yield strong oracle-type properties and, in certain settings (particularly orthonormal designs), behaves similarly to an $l_0$ penalty. Model selection methods such as the one discussed in [\citen{Robust_ident}] attempt to solve the non-relaxed, non-convex $l_0$ optimization problem using a subspace pursuit approach, though theoretical guarantees for such methods require assumptions that are often not met in practice. Prior work has addressed the challenge of applying SINDy on noisy data, including smoothing splines [\citen{spline_smooth_sindy_1}], interpolation [\citen{sindy_numerical_derv}], and integral formulations of sparse regression (e.g., weak SINDy [\citen{W_sindy,W_sindy_PDE}], D-SINDy [\citen{D_sindy}], etc.). 

Such sparse regression methods are effective for applications in which the differential variables are identified, all variables are measured, and the model is linear in the unknown parameters $\mathbf{p}$. However, for complex systems described by DAEs, one or more of these assumptions break down. Recent approaches have used nonlinear optimization for sparse-model selection in systems with partial measurements [\citen{DAHSI_paper,sindy_latent,stepaniants2024discovering,lu2022discovering}]. Yet, they rely on the assumption that the system is well-described by ordinary differential equations (ODEs) and do not explicitly incorporate algebraic constraints. Other groups have focused on model discovery for DAE systems but only when the conservation laws, network topology, or other algebraic constraints are known [\citen{dae_sindy_nandakumar,dae_power_grid_nandakumar,dae_symbolic_sari}].  Some DAEs can be converted to ODEs through substitution of the algebraic equations, reducing the dimension of state space and producing ODEs with rational functions on the right-hand side. Implicit-SINDy [\citen{I_sindy}], and the more computationally efficient parallel implicit SINDy (SINDy-PI) [\citen{PI_sindy}], have demonstrated sparse-model selection for ODEs with rational functions. However, these methods require that the differential variables are identified {\sl a priori} (see Section \ref{sec:challenges} for further discussion). A related concept was introduced in [\citen{sindy_feature_svd}], where principal component analysis (PCA) is used to eliminate correlated features from the library by retaining only those with the highest loadings in the first few principal components. However, this approach does not necessarily ensure a compact library, as the principal components may have uniformly distributed loadings unless a sparse PCA method is applied. No existing methods address the case where all states are measured, but the identity of the differential and algebraic variables is unknown, and the constraints must be identified.

In this paper, we propose Sparse Optimization for Differential-Algebraic Systems (SODAs), a novel approach for discovering DAE systems by sequentially identifying the algebraic equations and the differential equations. First, we iteratively identify algebraic relationships while partially parallelizing the optimization. Once these relationships are discovered, we refine the candidate function library to mitigate multicollinearity induced by algebraic relationships. Next, the dynamic component of the DAE system is discovered using the refined candidate function library and a selection of existing methods for ODE discovery ([\citen{sindy,W_sindy,IDENT_1,Robust_ident}]). This sequential approach improves the condition number of the library by eliminating redundant terms and restricts the challenge posed by noisy numerical derivatives in the ODE discovery phase. We employ a singular value decomposition (SVD) analysis of the candidate library as a guiding tool to determine the stopping criteria for the refinement process under high-noise conditions. The performance of our method is illustrated using numerous DAE systems, including chemical reaction networks (CRNs), multi-body mechanical systems, and power grids.

The paper is organized as follows. Section \ref{Sec:Background} provides a brief introduction to DAEs and the challenges associated with the existing SINDy framework for their discovery. The SODAs method is presented in Section \ref{Sec:Methods}, together with a detailed algorithm, and an example application to a simple CRN. In Section \ref{Sec:Results}, we demonstrate: 1) successful rediscovery of the DAEs used to simulate time series of chemical concentrations in CRNs, 2) rediscovery of the DAEs used in power grid modeling for 3 IEEE benchmark cases, and 3) discovery of a reduced coordinate system from pixel data obtained from real-time footage of single-pendulum experiments and from video recordings of animated single and double pendulums, which can potentially be leveraged for the subsequent discovery of ODEs. Finally, in Section \ref{Sec:DiscussionSec}, we address various challenges associated with our method and discuss possible improvements.

\section{Background: DAEs and existing SINDy frameworks}\label{Sec:Background}

\subsection{Differential-Algebraic Equations}\label{SEc:DAE}
In the types of physical systems in which we are interested, DAEs commonly arise from quasi-steady state approximations, conservation laws, and physical constraints, as seen in the modeling of chemical reactions, power grids, and multi-body systems [\citen{DAE_applications}]. When DAEs arise from quasi-steady state approximations [\citen{quasi_steady_chem,quasi_steady_diffusion}], where fast variables are assumed to reach equilibrium instantaneously relative to other variables operating in slower dynamics, it is often challenging to identify which variables are in quasi-steady state and are, therefore, algebraic rather than differential. 

Let $\mathbf{x}(t)=(\mathbf{y}(t),\mathbf{z}(t))$ represent the state variables in the dynamical system, such that $\mathbf{y}=[y_1,\ldots,y_m]^T\in C^{c}([0,T],\mathbb{R}^{m})$ and $\mathbf{z}=[z_1,\ldots,z_{d-m}]^T\in C^{0}([0,T],\mathbb{R}^{d-m})$. The following set of equations represents a system of DAEs in a semi-explicit form, where the algebraic constraints can be separated from the differential operator:  
\begin{align}
\mathbf{F}(\mathbf{y},\dot{\mathbf{y}},\ddot{\mathbf{y}},\ldots,\mathbf{z},t) & = 0, \label{eq:DAE_main_1}\\
\mathbf{G}(\mathbf{y},\mathbf{z},t) & = 0, \label{eq:DAE_main_2}
\end{align}
where $\mathbf{F}=(f_{1},\ldots,f_{m})$ and $\mathbf{G}=(g_{1},\ldots,g_{d-m})$ represent the differential and algebraic parts of the system, respectively.  Note that time derivatives, up to order $c$, appear only in the differential state variables $\mathbf{y}$, while the algebraic variables $\mathbf{z}$ are continuous on $[0,T]$. DAEs can also be expressed in an implicit form, where the algebraic constraints are not separated, and our algorithm is designed to work for both forms. Note that the implicit form can be converted into a semi-explicit form and also that the latter form is often more natural for many applications, as it reflects the physical significance of the explicit algebraic constraints present in these contexts. All DAE examples discussed here and used to test the algorithm have a semi-explicit form.

The DAE system (\ref{eq:DAE_main_1})--(\ref{eq:DAE_main_2}) is well-posed depending on certain regularity conditions on $\mathbf{F}$ and $\mathbf{G}$ as well as the existence of compatible initial conditions $\mathbf{x}_0$ satisfying the algebraic constraints (\ref{eq:DAE_main_2}) (see [\citenum{DAE_book_1}], Chapter 2). DAEs can be further classified based on \textit{indices} related to the easiness of conversion to an ODE system (via subsequent differentiation) or based on the complexity of the numerical solution. For a detailed analysis of the solvability of nonlinear DAE systems and their classification based on various indices, see [\citenum{DAE_book_1}] and [\citenum{DAE_numerical_2}]. We assume that the DAE system to be discovered from data is well-posed and that the state-space trajectories $\mathbf{x}(t)$ are smooth enough on the interval $[0,T]$. 

A simple example of DAE system describes the chemical kinetics of an enzyme-mediated, irreversible reaction:
\begin{align}
&\frac{d[A]}{dt}=  -k_{1}[E_1][A]+k_{2}[AE_1], \label{eq:CRN-sim-1}\\
&\frac{d[B]}{dt}=  k_3[AE_1], \label{eq:CRN-sim-2} \\
&[E_1]+[AE_1]-E_1^{\text{tot}}= 0, \label{eq:CRN-sim-3} \\
&\frac{d[AE_1]}{dt} = -(k_2+k_3)[ES]+k_1[E][S]\approx  0, \label{eq:CRN-sim-4}
\end{align}
where $[A], [E_1], [AE_1], [B]$ are concentrations of substrate, enzyme, enzyme-substrate complex, and products of a CRN, respectively. The enzyme-substrate complex $[AE_1]$ can be assumed to be in quasi-steady state under many experimental conditions [\citen{segel1989quasi}]. Once \eqref{eq:CRN-sim-4} is assumed to be in quasi-steady state, the differential variables are $[A]$ and $[B]$, and the algebraic variables are $[E_1]$ and $[AE_1]$. 

\subsection{Data-driven discovery of equivalent ODE systems}\label{sec:challenges}

An approach previously explored is to discover the DAEs in their equivalent reduced ODE form. DAE systems of the form~\eqref{eq:DAE_main_1}--\eqref{eq:DAE_main_2} can often be converted to a fully coupled ODE system:
\begin{align}
\widehat{\mathbf{F}}(\mathbf{y},\dot{\mathbf{y}},\ddot{\mathbf{y}},\ldots,\mathbf{z},t) & = 0. \label{eq:DAE_reduced}
\end{align}

\noindent
This reduction is performed regularly in many fields, including chemical, mechanical, and electrical engineering [\citen{DAE_applications,quasi_steady_chem}], to reformulate models in a manner compatible with well-established techniques for ODEs, such as numerical solvers and control methods. The reduction involves eliminating state variables using algebraic constraints and can lead to more complicated expressions in~\eqref{eq:DAE_reduced} as compared to~\eqref{eq:DAE_main_1}--\eqref{eq:DAE_main_2}. For example, the system~\eqref{eq:CRN-sim-1}--\eqref{eq:CRN-sim-4} can be reduced using \eqref{eq:CRN-sim-3} and \eqref{eq:CRN-sim-4} to eliminate $[E_1]$ and $[AE_1]$ resulting in
\begin{align}
\frac{d[A]}{dt} &= -\frac{k_1  k_3  E_1^{\text{tot}}  [A]}{k_2 + k_3 + k_1  [A]}. \label{eq:CRN-sim-redODE}
\end{align}
These are the canonical Michealis-Menten kinetics [\citen{johnson2011original,srinivasan2022guide}]. More generally, DAE systems of monomial terms can often be reduced to the following form, especially when the algebraic variables can be explicitly solved from the algebraic equations 
\begin{align}
\frac{dx_i}{dt} &=\frac{\sum_{j_n} p_{j_n}^i n_{j_n}^i(\mathbf{x})}{\sum_{j_d}p_{j_d}d_{j_d}^i(\mathbf{x})}, \quad i = 1, \ldots, d, \label{eq:DAE-red-example}
\end{align}
where $n_{j_n}^i(\mathbf{x}),  d_{j_d}^i(\mathbf{x}) \in \mathbb{R}[\mathbf{x}]$ are multivariate monomials of the state vector $\mathbf{x}$. This reduction to ODEs has resulted in rational functions of the state variables compared to the original set of equations, which only contained polynomial functions. Converting higher-index DAEs to a reduced ODEs can become more complicated, resulting in rational functions with nonlinear terms in both the numerator and denominator, even when the original DAE system consists of simple additive terms. 

Although system \eqref{eq:DAE-red-example} is now converted to an ODE form, the optimization cost \eqref{eq:loss-cont} formulated in SINDy is not a convex problem anymore because the functions are rational and hence the model is not linear in the coefficients $\textbf{p}$. This nonlinear optimization can be challenging due to various factors including non-convexity, poles of rational functions, and stability issues for techniques that use an ODE integration as part of the loss function evaluation. 
To address this issue, implicit-SINDy ~[\citen{I_sindy}] reformulated the reduced ODEs of the form~\eqref{eq:DAE-red-example} as an implicit ODE and then applied sparse optimization. Expressed implicitly,~\eqref{eq:DAE-red-example} becomes
\begin{align}
\frac{dx_i}{dt}\sum_{j}d_j^i(\mathbf{x}) - \sum_{j}n_j^i(\mathbf{x}) &=0, \quad i = 1, \ldots, d. \label{eq:Impl-DAE}
\end{align}
In this implicit-SINDy framework, the candidate function library is augmented to be $\tilde{\candLibMatrix}=\{\dot{\mathbf{x}}\candLibMatrix,\candLibMatrix\}$, where $\dot{\mathbf{x}}\candLibMatrix$ denotes the set of all products of the state derivatives with the candidate functions. As a result, the library has $2J$ terms. The goal is then to satisfy equation~\eqref{eq:Impl-DAE} by finding a sparse vector in the nullspace of this augmented library $\tilde{\candLibMatrix}$, which can be achieved via the alternating direction method of multipliers (ADMM). Nonetheless, this method has some drawbacks, such as its sensitivity to noise when estimating the nullspace and its high data requirements.

The SINDy-PI method~[\citen{PI_sindy}] improves upon the implicit-SINDy method by performing a partial combinatorial search: each library element in $\tilde{\candLibMatrix}$ is moved to the right-hand side of the ODE and standard sparse regression is applied assuming the moved term is part of the true equation. Solving the $2J$ convex, sparse linear regressions can also be parallelized in this approach. Notwithstanding, these implicit methods still require large amounts of data, and the robustness to noise is often insufficient for many practical applications. Specifically, SINDy-PI was applied to the Michaelis-Menten kinetics model using a large dataset consisting of 150,000 data points for training and testing, generated from 3,000 randomly sampled initial conditions. Despite this large volume of data, the method was shown to handle only Gaussian noise up to a magnitude of $\sigma=0.04$ (Figure 3 in [\citen{PI_sindy}]), which corresponds to a noise level of less than 1\% given that the initial condition for the chemical concentration ranged between 0 and 12.5 in magnitude. 

\subsection{Benefits and challenges of model discovery in the DAE framework}\label{Sec:benefitsDAESec}

In addition to the high data quality requirements, a fundamental issue is that implicit-SINDy and SINDy-PI methods do not identify the algebraic variables and associated equations within a DAE system. The algebraic constraints in the DAEs typically carry physical significance related to physical constraints or quasi-steady-state approximations, and this intuitive meaning is lost in the equivalent ODE formulation. Creating an explicit system of ODEs in the form \eqref{eq:DAE_reduced} can be complex, computationally intensive, or even impossible for higher-index DAE systems or when the algebraic constraints are not differentiable in time [\citen{petzold_paper}]. Furthermore, small parameter changes in the original DAE system often lead to distinct reduced ODE forms. Solving or discovering DAEs directly preserves the underlying structure of the problem. 
An additional benefit of identifying the algebraic relationships is that the library can be iteratively simplified by removing an algebraic variable for every algebraic equation discovered. Removal of terms with algebraic variables reduce multicollinearity among the candidate library terms due to the explicit algebraic constraints, improving the conditioning of the library. Consequently, the discovery of DAEs provides a more interpretable and computationally beneficial framing.

Discovering a DAE system in its semi-explicit form 
(\ref{eq:DAE_main_1})--(\ref{eq:DAE_main_2}) also offers flexibility for numerical analysis and simulation. For example, in CRN applications, simulations are often carried out using ODE solvers designed for stiff systems; however, in applications involving high-index DAEs, solving a reduced ODE system can lead to numerical instability [\citen{petzold_paper,DAE_kunkel, DAE_book_1}]. In these cases, single-step and multi-step methods for solving DAEs in their non-reduced forms are preferred [\citen{petzold_paper,DAE_kunkel, DAE_book_1, DAE_numerical_1, DAE_numerical_2}]. By first discovering the DAE system in the semi-explicit form, the user can choose either to convert it to a reduced ODE or to keep the original 
semi-explicit form, depending on which numerical solvers are most appropriate. In comparison, converting a reduced ODE system back into a DAE system with 
explicit constraints is often not feasible.


To address all the issues discussed so far, including the interpretability of physical constraints, the numerical instabilities caused by perfect multicollinearity in the candidate library, and the adversial effects of noise present in implicit ODE discovery, we propose the method SODAs. By separating the discovery of differential and algebraic equations, we directly identify the algebraic variables and equations and eliminate the impact of noisy derivatives on the discovery of algebraic equations. We can then utilize the knowledge gained from the discovered algebraic equations to better condition the library before attempting to discover differential equations.

\section{Methods: Sparse Optimization for Differential-Algebraic Systems}\label{Sec:Methods}

\begin{figure}[h!]
    \centering
    \includegraphics[width=1\linewidth]{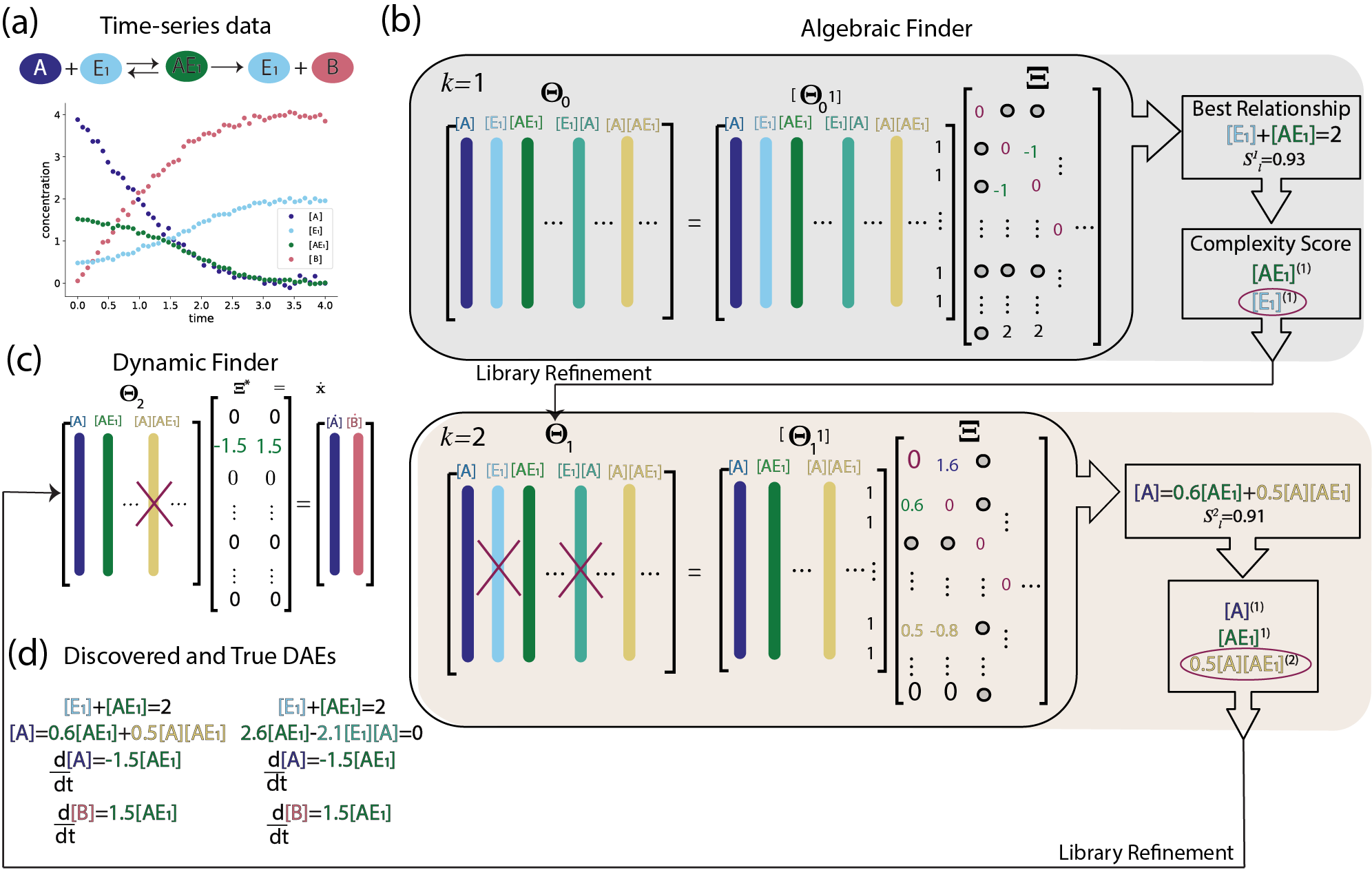}
    \caption{\textbf{SODAs algorithm explained using a chemical reaction network example.} 
    (\textbf{a}) Collect time series and construct candidate library matrix. 
    (\textbf{b}) Algebraic Relation Finder: iteratively finds algebraic relationships and refines the candidate library based on the discovered relationship. In each iteration $k$, the best relationship is selected based on its score $S^k$, and a complexity score is assigned to each term within this relationship. This complexity score is used to select a term to refine the candidate library for the next iteration: the red X-s indicate the removal of features from the candidate library that are multiples of the selected term.  
    (\textbf{c}) Dynamic Finder: the refined candidate library from the algebraic finder step is used to find the system of ODEs.
      (\textbf{d}) Assembled DAEs: The algebraic and differential equations from the algebraic and dynamic finder steps respectively, are assembled to form the final set of DAEs (\textit{left}), which is equivalent to the true set of DAEs (\textit{right}).}
    \label{fig:sodas_main_algo}
\end{figure}

In this section, we introduce our algorithm SODAs (Figure \ref{fig:sodas_main_algo}, Algorithm 1), and demonstrate how it can be used to discover DAEs in their explicit form~\eqref{eq:DAE_main_1}--\eqref{eq:DAE_main_2}. Details on the inputs, parameters, and defaults used in Algorithm 1 can be found in Supplementary-section S1. The core idea is to separate the discovery of the algebraic and dynamic state equations. \href{https://pypi.org/project/DaeFinder}{\texttt{DaeFinder}} is a Python implementation of SODAs that allows users to easily create popular candidate libraries relevant to many applications and efficiently implement various components of the SODAs algorithm.  Details about this package and defaults used in the algorithm can be found in Supplementary-section S2. SODAs iteratively identifies algebraic relationships, removes them from the library, and reduces multicollinearity within the library. After this process, SINDy is applied to discover the dynamic equations from the remaining library terms. We illustrate the application of SODAs using an example of a simple CRN~\eqref{eq:CRN-sim-1}--\eqref{eq:CRN-sim-4} (Figure~\ref{fig:sodas_main_algo}a). 

\subsubsection*{Method formulation}

Let \(\mathbf{X} = [\mathbf{X}_1 \; \mathbf{X}_2 \; \ldots \; \mathbf{X}_d]\) denote the \(N \times d\) time-series data matrix, where each column \(\mathbf{X}_i\) contains possibly noisy measurements of the state \(x_i\) at time points \(\{t_1 = 0, \ldots, t_N = T\}\), such that \(\mathbf{X}_{ij} \approx x_j(t_i)\). Note that the measurements \(\textbf{X}_{ij}\) may also be obtained through transformations of the original state variables. However, for simplicity, throughout this paper, we assume that all states are directly observable with some measurement noise. 
SODAs is typically applied to DAE systems where the terms appearing in equations \eqref{eq:DAE_main_1}--\eqref{eq:DAE_main_2} can be written as linear combination of a set of predefined library terms which comes from prior knowledge. We start by defining a candidate library  \(\candLibMatrix = \{ \theta_1(\mathbf{x}), \ldots, \theta_J(\mathbf{x}) \}\), which is a set of functions of the state variables. 
Based on the time-series data $\mathbf{X}$, we next form the initial candidate library matrix denoted by  \(\candLibMatrix_0= [ \theta_1(\mathbf{X})\ \theta_2(\mathbf{X})\ \ldots\ \theta_J(\mathbf{X}) ]\), which is an $N\times J$ matrix where the $j^{\text{th}}$ column represents the evaluation of the candidate function $\theta_j$ at time points \(\{t_1 = 0, \ldots, t_N = T\}\). In our applications, a polynomial library is used for the CRN and the nonlinear pendulum, while a combination of polynomial and trigonometric functions is used for the power grid. After forming the candidate library, SODAs consists of two major steps: the algebraic finder step (Figure~\ref{fig:sodas_main_algo}b) and the dynamical system finder step (Figure~\ref{fig:sodas_main_algo}c).

\subsubsection*{Identifying algebraic relationships through iterative sparse optimization}

 For many systems, it is unclear {\sl a priori} which states are differential variables and which are algebraic, especially if separation of timescales is expected and a subset of variables are in quasi-steady state. Our approach does not necessarily require pre-existing knowledge about which states are in the algebraic relations. Notwithstanding, prior knowledge on state variables that take part in algebraic equations can be incorporated into the pipeline. For unknown algebraic equations, the algebraic finder step involves iteratively identifying algebraic relationships and reducing multicollinearity in the library by removing these relationships. Once all algebraic equations are found,  the remaining dynamical equations are identified from the refined library using an existing method like SINDy [\citen{sindy,D_sindy,W_sindy}] or IDENT [\citen{IDENT_1,Robust_ident}]. 



At each $k^{\text{th}}$ iteration of the algebraic finding step, we start with a refined candidate library matrix $\candLibMatrix_{k-1}$, where $\candLibMatrix_0$ represents the initial candidate library matrix $\candLibMatrix$ (Figure \ref{fig:sodas_main_algo}b). Let $\mathcal{I}_k$ represent the index set for elements of $\candLibMatrix_k$ and $J_k = |\mathcal{I}_k|$.
Assuming no prior information about the number of algebraic relationships and which state variable functions are part of algebraic relationships, we sparsely fit each member of the candidate library $\candLibMatrix_{k-1}$ against the rest of the library. That is,
\begin{align}
\theta_l(\mathbf{x}) &= \sum_{\substack{j\in \mathcal{I}_k\\ j\ne l}} p_j^l \theta_j(\mathbf{x}), \quad l\in  \mathcal{I}_{k-1}. \label{eq:sodas-algo-1}
\end{align}
To enforce sparsity, such that most of the regression parameters \(p_j^l = 0\) for \(j \in \mathcal{I}_{k-1}\), the following loss function can be minimized:
\begin{align}
L_l(\mathbf{p}^l) &=\sum_{r=1}^{N} \Big( \theta_l(\mathbf{x}(t_r)) - \sum_{\substack{j\in \mathcal{I}_k\\ j\ne l}} p_j^l \theta_j(\mathbf{x}(t_r)) \Big )^2 + \alpha \|\mathbf{p}^l\|_0, \quad l\in  \mathcal{I}_{k-1}. \label{eq:sodas-algo-2}
\end{align}
In matrix form, this series of regression problems can be written as
\begin{equation}
\candLibMatrix_{k-1} = \candLibMatrix_{k-1}\Xi, \label{eq:sodas-algo-matrix-form}
\end{equation}
where $\Xi = \left[\mathbf{p}^1\ \mathbf{p}^2\ \ldots\ \mathbf{p}^{J_k} \right]$ with constraints $\Xi_{jj}=0$ and $\|\Xi\|_0$ is minimized (Figure \ref{fig:sodas_main_algo}b, \textit{ left}). Note that, in practice, the matrix form is not used, but instead each $\mathbf{p}^l$ is optimized in parallel by minimizing the loss function $L_l(\mathbf{p}^l)$. A similar parallelized approach is introduced in~[\citen{PI_sindy}] for semi-combinatorial discovery of implicit equations of the form~\eqref{eq:Impl-DAE}, whereas in SODAs this parallel fitting approach is used for finding algebraic relationships. 
Notably, derivatives are not involved in the algebraic finder step and do not need to be estimated from the noisy data.
There are multiple established approaches to enforce the sparsity of the parameter matrix $\Xi$.  The hyperparameters used in this step include the thresholding parameter $\alpha_{\text{th}}$ [\citen{sindy}], which determines the coefficient size for removing terms, as well as the regularization penalty weight $\alpha$ [\citen{Lasso_tibshirani}], both of which  are encoded in $\alpha^a$. In practice, these hyperparameters should be chosen using a Pareto-front  approach, as discussed in [\citen{sindy}].

\subsubsection*{Improving numerical conditioning of the candidate library}

The next step of the algorithm is to select the best-fitting algebraic relationship across all the candidate fits in \eqref{eq:sodas-algo-2} (Figure~\ref{fig:sodas_main_algo}b, \textit{right}). To identify the most likely relationship, we assign a score $S^k_l$ to each candidate fit in \eqref{eq:sodas-algo-2}. In our examples, we have used the coefficient of determination $\mathcal{R}^2$  [\citen{ISL_book}], which for a given feature  $\theta_l$ with mean $\bar{\theta_l}$  is given by  $\mathcal{R}^2=1-\frac{ \sum_{r=1}^{N} \Big ( \theta_l(\mathbf{x}(t_r)) - \sum_{\substack{j\in \mathcal{I}_k\\ j\ne l}} p_j^l \theta_j(\mathbf{x}(t_r)) \Big )^2 }{ \sum_{r=1}^{N} \Big ( \theta_l(\mathbf{x}(t_r)) - \bar{\theta_l} \Big )^2}$, and ranges from 0 to 1, where larger values indicate a better fit. Other criteria, such as AIC and BIC (see~[\citen{mangan_IC,aic_bic_chakrabarti}]), that balance training error with model simplicity can also be used to define $S^k_l$. The fit with the highest score is then identified as the most prominent algebraic relationship in the library $\candLibMatrix_{k-1}$ and is further used to refine $\candLibMatrix_{k-1}$. In our experience, conservation laws that involve first-order terms are initially identified as the best relationships during this process.

After identifying the best algebraic relationship based on the score \(S^k_l\) for the \(k^\text{th}\) iteration, any common factors in that relationship are factored out to obtain a reduced equation \(g_k(\mathbf{x}) = 0\). To address the multicollinearity introduced by this algebraic relationship, we remove one of the terms in \(g_k(\mathbf{x})\). To decide which term to remove, each term \(\theta^{(k)}\) in \(g_k(\mathbf{x})\) is assigned a complexity score (Figure~\ref{fig:sodas_main_algo}b, \textit{right}). In polynomial libraries, this score is set to be the degree of the term by default, i.e., higher-degree monomials should have a higher complexity score. The term with the highest complexity score,  \(\theta^{(k)}_{\text{high }}\) (circled in Figure~\ref{fig:sodas_main_algo}b, \textit{right}), is then removed from the library. The rationale behind this refinement is that if \(\theta^{(k)}_{\text{high }}\) is part of the discovered algebraic relationship \(g_k(\mathbf{x}) = 0\), then \(\theta^{(k)}_{\text{high }}\) can be expressed as a linear combination of the other terms in \(g_k(\mathbf{x})\), so there is no loss of information by removing a term from the library. By keeping lower rather than higher complexity terms in the library, we reduce noise amplification [\citen{messenger2022asymptotic}].  If multiple terms share the same complexity score, then either expert knowledge or quality of measurement of states (through statistical measures like variance) can also be used to choose one term over the other, or a term can be selected randomly. Additionally, for any \(\theta \in \candLibMatrix_{k-1}\), the equation \(\theta \, g_k(\mathbf{x}) = 0\) also represents a potential algebraic relationship that can lead to perfect multicollinearity. To remove all degenerate relationships of the form \(\theta \, g_k(\mathbf{x}) = 0\), all multiples of \(\theta^{(k)}_{\text{high }}\)are also removed from \(\candLibMatrix_{k-1}\).
Let $\Gamma_k$ be a subset of the candidate library matrix $\candLibMatrix_{k-1}$ that has \(\theta^{(k)}_{\text{high }}\) as a factor; then $\candLibMatrix_{k} = \candLibMatrix_{k-1} \setminus \Gamma_k$, where all elements of $\Gamma_k$ are removed from $\candLibMatrix_{k-1}$ as part of the refinement process. It follows that
the libraries at each refinement step form a decreasing chain, $\candLibMatrix_{k}\subset\candLibMatrix_{k-1}\subset\cdots\subset\candLibMatrix_{0}.$

\subsubsection*{Mitigating multicollinearity due to algebraic relationships }\label{sec:stopping_crit}

As an example, in the CRNs \eqref{eq:CRN-sim-1}--\eqref{eq:CRN-sim-4}, the algebraic relationship $[E_1]+[AE_1]-E_1^{tot}= 0$ can cause the following perfect multicollinearity in the candidate library matrix $\candLibMatrix$ that has up to degree 2 monomials:
\begin{align}
&[E_1]+[AE_1]-E_1^{tot}= 0, \label{eq:CRN-alg-1}\\
&[E_1]\left([E_1]+[AE_1]-E_1^{tot}\right)= 0, \label{eq:CRN-alg-2}\\
&[A]\left([E_1]+[AE_1]-E_1^{tot}\right)= 0, \label{eq:CRN-alg-3}\\
&[B]\left([E_1]+[AE_1]-E_1^{tot}\right)= 0, \label{eq:CRN-alg-4}\\
&[AE_1]\left([E_1]+[AE_1]-E_1^{tot}\right)= 0. \label{eq:CRN-alg-5}
\end{align}
For higher-order candidate libraries, there will be more relationships. All such degenerate algebraic relationships make up the nullspace of the library matrix $\candLibMatrix$. If any one of the degenerate equations~\eqref{eq:CRN-alg-2}--\eqref{eq:CRN-alg-5} is identified as the relationship with the best $S_l^k$, it is first reduced to the form~\eqref{eq:CRN-alg-1}. Both terms in this reduced relationship (i.e., $[E_1]$ and  $[AE_1]$) have the same complexity score of 1 based on the degree of the monomial (Figure~\ref{fig:sodas_main_algo}b, \textit{top-right}). In the absence of any other information, we randomly choose to select $[E_1]$ to refine the library. Consequently, all factors of $[E_1]$ including $[E_1]$, $[A][E_1]$, $[B][E_1]$, $[AE_1][E_1]$, and $[E_1]^2$ are removed from the candidate library $\candLibMatrix_0$  (see how $\candLibMatrix_1$ is defined for $k=2$ in Figure~\ref{fig:sodas_main_algo}b). 
This process of algebraic equation discovery and library refinement are iteratively continued until all algebraic relationships are discovered.  Note that the algebraic equations found will not be unique because any linear combination of algebraic equations can lead to another algebraic relationship, but our approach tries to find a minimal set of equivalent algebraic relationships that can form any algebraic relationship using linear combinations.

\begin{figure}[h!]
    \centering
    \includegraphics[width=0.85\linewidth]{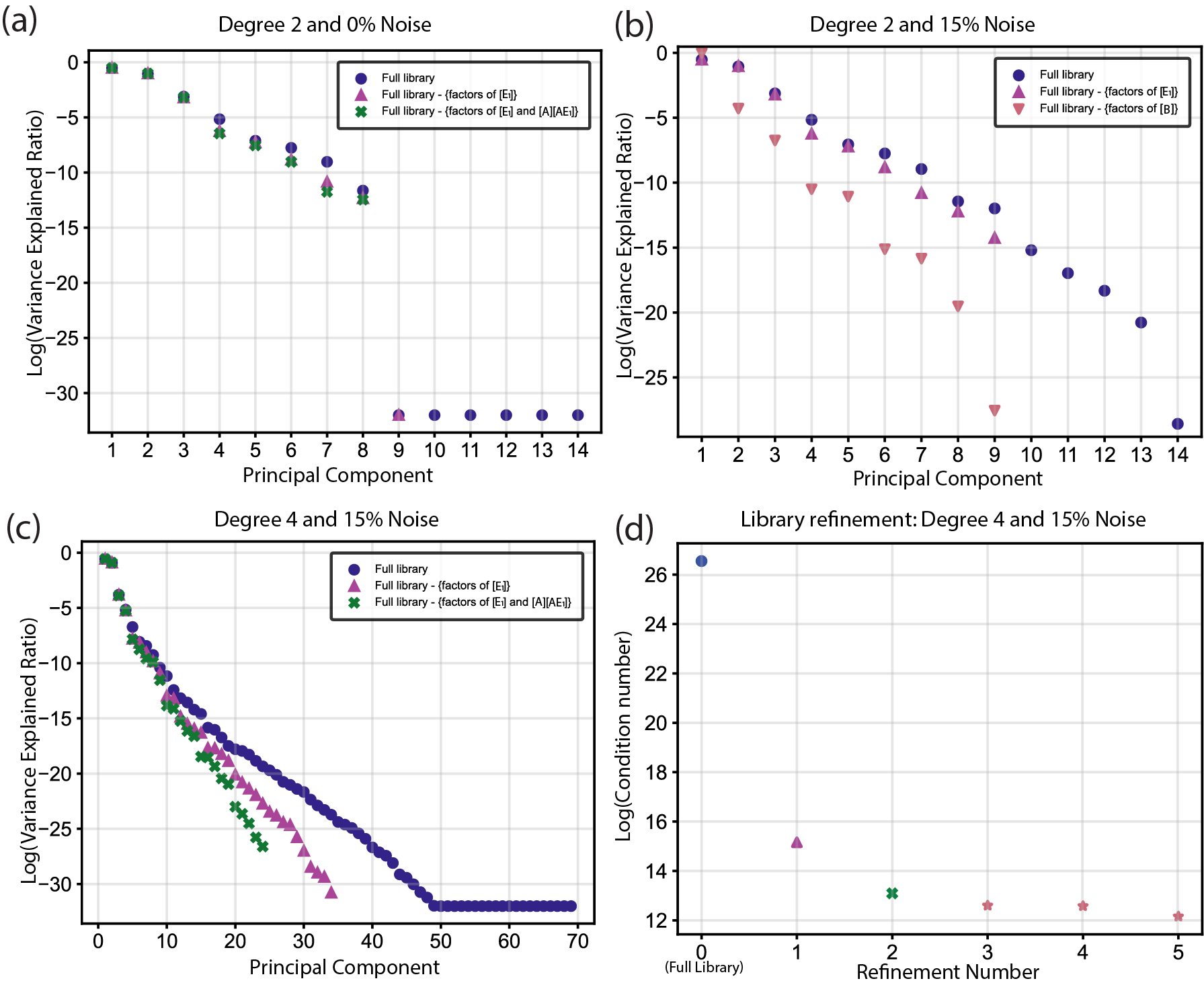}
    \caption{\textbf{SVD analysis of the candidate library \(\candLibMatrix\)}: The Y-axes are presented on a natural logarithm scale. (\textbf{a}) library with monomials up to degree 2 and no noise, (\textbf{b}) library with monomials up to degree 2 and 15\% noise, (\textbf{c}) library with monomials up to degree 4 and 15\% noise, (\textbf{d}) tracking condition number after each refinement for library with monomials up to degree 4 and 15\% noise.  SVD analysis helps track redundancies in the candidate library space and determines when to stop iterations in the algebraic finder step.}
    \label{fig:svd-analysis}
\end{figure}

For many applications, such as power-grid networks, the number of algebraic constraints in a DAE system is known beforehand from physical principles such as the conservation of energy at each node. However, for other applications, such as CRNs, the number of algebraic constraints is not known in advance. When the number is unknown, an SVD analysis of the candidate library matrix \(\candLibMatrix_{k-1}\) can be used to track the number of algebraic relationships remaining in the library (Figure~\ref{fig:svd-analysis}) at iteration $k$. Consider the SVD decomposition 
\begin{equation}
\candLibMatrix_{k-1} = \mathbf{U}\boldsymbol{\Sigma}\mathbf{V}^\top,
\end{equation}
where \(\boldsymbol{\Sigma} = \mathrm{diag}(\sigma_1,\ldots,\sigma_r)\) contains the singular values \(\sigma_1 \ge \cdots \ge \sigma_r \ge 0\), with \(r = \mathrm{rank}(\candLibMatrix_{k-1})\). The right singular vectors associated with zero singular values ($\sigma_i \approx 0$) span the nullspace of \(\candLibMatrix_{k-1}\) and correspond to exact linear dependencies in the library. Consequently, the nullity
 $\mathrm{nullity}(\candLibMatrix_{k-1}) = J - r$
provides an estimate of the number of linear algebraic constraints that exist in the library. To quantify the relative contribution of each principal direction, we calculate the variance-explained ratio associated with the \(i\)-th singular value using 
\begin{equation}
\eta_i = \frac{\sigma_i^2}{\sum_{j=1}^{r} \sigma_j^2}. 
\end{equation}
When noise levels are low, the spectrum \(\{\eta_i\}\) shows a distinct separation between the dominant principal components and directions in the null space. As valid algebraic relationships are iteratively removed from the library, the zero singular values disappear, while the number of dominant components stays the same until all algebraic constraints have been eliminated. (Figure~\ref{fig:svd-analysis}a).

Even at high noise, where the separation becomes less clear, the SVD spectrum can indicate whether the refinement is based on a true algebraic relationship in the library. For instance, at 15\% noise (see \eqref{eq:noise} to see how noise is added to time series) in Figure~\ref{fig:svd-analysis}b, refining the library using the term \([E]\) (which is part of an algebraic relationship) reduces the number of zero singular values while keeping the number of dominant principal directions, those with relatively large variance explained,  unchanged. In contrast, refining it using \([P]\) (which is not part of any algebraic relationship) causes a loss in the number of dominant principal directions. Through a correct refinement, the candidate library can become significantly smaller and better-conditioned, aiding the discovery of dynamic equations. The effect of a series of correct library refinements in removing only low-variance-explained-ratio values is more pronounced for libraries with higher-order terms (e.g., monomials up to degree 4, see Figure~\ref{fig:svd-analysis}c).
SODAs can be used with or without prior knowledge of the number of algebraic relationships, $K$. If $K$ in the DAE system is not known beforehand, it can determined via a Pareto front analysis of the library's condition number improvements. We use the lack of improvement in the condition number within some tolerance as a stopping criterion for the algebraic discovery and refinement process ie, $\text{Cond}(\candLibMatrix_{K+1})\approx\text{Cond}(\candLibMatrix_{K})<\text{Cond}(\candLibMatrix_{K-1})<\cdots<\text{Cond}(\candLibMatrix_0)$. For example, the improvement in conditioning stagnates after the first two refinements for the CRN example, and so the Pareto front suggests stopping library refinement after the second equation is removed (Figure~\ref{fig:svd-analysis}d). For very noisy data it may be necessary to treat K as a hyperparameter. If some algebraic relationships are known ahead of time, then SODAs ensures at least that many relationships are recovered.

\subsubsection*{Identifying differential variables and their differential equations}

Once the algebraic relations are found, the differential equations remain to be discovered. However, we must first assign which state variables will be treated as differential variables. Each algebraic relationship will impose algebraic constraints on one or more of the state variables. If multiple variables are part of an algebraic relationship, domain expertise or quality of measurements can be used to decide which variables should be identified as the differential variables in the DAEs. For example, the concentration of enzyme complexes in a CRN are often selected as the algebraic variables because they are difficult to measure directly, and it is advantageous to reduce them out of the system through substitution of the algebraic equations as in \eqref{eq:CRN-sim-redODE}. Therefore, metabolite concentrations may be preferred as the differential variables. Another consideration is the quality of measurements and the ability to find accurate time derivatives of state variables. If a state variable exhibits a relatively high dynamic range in measurement and we have access to high-quality, high-sampling measurements of this state, then this state can be treated as a differential variable.

For the identified differential variables, any existing methods for ODE discovery (e.g., SINDy) can be used to find the dynamics equations in the DAE system (Figure \ref{fig:sodas_main_algo}c). This is the stage where noisy derivatives affect recovery, and there are various approaches to mitigate noise, including weak formulation approaches [\citen{W_sindy,W_sindy_PDE,D_sindy}] and smoothing functions  [\citen{spline_smooth_sindy_1,savitzky_filter}]. We have used Savitzky-Golay filters [\citen{savitzky_filter}] in our examples. For numerical experiments, we have used LASSO followed by sequential thresholding to solve the sparse optimization. 
Finally, we obtain a set of differential and algebraic equations of the form~\eqref{eq:DAE_main_1}--\eqref{eq:DAE_main_2} (Figure \ref{fig:sodas_main_algo}d).

\begin{algorithm} \label{alg:sodas}
\caption{SODAs: Sparse Optimization for Differential-Algebraic Systems}
\KwIn{$\mathbf{X}$,  $K$, $\candLibMatrix$, $\alpha^a$, $\alpha^d$, \texttt{sparse\_fit}, \texttt{score\_function}, $\epsilon$}
\KwOut{ discovered DAEs
}
\BlankLine

\textbf{Step 1: Construct library matrix:} \\
$\candLibMatrix_0 = [\theta_1(\mathbf{X}), \ldots, \theta_J(\mathbf{X})]$ \\
\BlankLine

 \textbf{Step 2: Algebraic Finder} \\
$k=0$\quad \quad \textit{\%set iteration index}\\
\texttt{condition\_number\_improved = True} \quad \quad \textit{\%flag for condition number improvement. }\\

\While{($k\le K$ or \texttt{condition\_number\_improve}d)}{
   $k=k+1$\\
    \ForEach{$\theta_l \in \candLibMatrix_{k-1}$ \textit{\quad (done in parallel)}}{
       \texttt{sparse\_fit}$(\candLibMatrix_{k-1},\theta_l, \alpha^a)$ \quad \quad \textit{\%sparse  regression  $\theta_l(\mathbf{x}) = \sum_{\substack{j \in \mathcal{I}_k \\ j \neq l}} p_j^l \theta_j(\mathbf{x}) +\alpha \|\mathbf{p}^l\|_0$}.\\
    $S^k_l$ = score\_function$(\mathbf{p}^l,\candLibMatrix_{k-1})$
    }

    \textbf{Library Refinement Steps:}\\
    $m = \arg\max\limits_{l}(S^k_l)$\quad \quad \textit{\%find the best relationship based on the maximum score}.\\
    
    $g_k(\mathbf{x}) = \texttt{equation\_reduce}(m, \candLibMatrix_{k-1})$\quad \quad \textit{\%reduce the equation $\theta_m(\mathbf{x}) = \sum_{j} p_j^m \theta_j(\mathbf{x})$ to $g_k(\mathbf{x}) = 0$.}\\
  $\theta^{(k)}_{\text{high}} = \texttt{select\_term}(g_k(\mathbf{x}))$\quad \quad \textit{\%select the term in $g_k(\mathbf{x})$ based on complexity rank}.\\
    $\Gamma_k = \texttt{multiples}(\theta^{(k)}_{\text{high}},\candLibMatrix_{k-1})$ \quad \quad \textit{\%identify all multiples of $\theta^{(k)}$ in $\candLibMatrix_{k-1}$}.\\
  $\candLibMatrix_{k} = \candLibMatrix_{k-1} \setminus \Gamma_k$ \quad \quad \textit{\%update the library}.\\
  \texttt{condition\_number\_improved} = \texttt{bool}($\text{Cond}(\candLibMatrix_{k-1})-\text{Cond}(\candLibMatrix_{k})>\epsilon[k]$)\\
}
$K=k-1$ \\
$\candLibMatrix_*=\candLibMatrix_{k-1}$ \quad \textit{\%setting final library to the one }.\\
\BlankLine
\textbf{Step 3: Dynamical System Finder} \\
  $\mathbf{y} = \texttt{find\_dynamic\_states}(\mathbf{x}, \{g_k\}_{k=1}^K)$\quad \quad  \textit{\%finding dynamic states based on discovered algebraic equations.} \\
$ \mathbf{F}(\dot{\mathbf{y}},\ddot{\mathbf{y}},\mathbf{x},\ldots)=0\xleftarrow{}\texttt{discover\_ODEs}(\mathbf{y},\Theta_*, \alpha^d)$ \quad \quad  \textit{\%discovering ODEs  using refined library.} \\
\BlankLine

\Return $\{ \{g_k\}_{k=1}^K),\mathbf{F}\}$\quad \quad  \textit{\%return discovered DAE system.} \\
\begin{footnotesize}
Details on the inputs and defaults can be found in Supplementary-section S1.
\end{footnotesize}
\end{algorithm}

\section{Results}\label{Sec:Results}
We demonstrate the capability of SODAs to discover DAEs and identify reduced coordinate systems through three distinct examples. In Example 1, SODAs rediscovers the underlying DAEs of CRNs with varying complexity, demonstrating the identification of conservation laws and quasi-steady states. In Example 2, SODAs rediscovers DAEs used in power-grid modeling for three benchmarks (IEEE-4, IEEE-9, and IEEE-39) demonstrating recovery of network topology when the library contains trigonometric functions.  Example 3  applies SODAs to pixel data from real-time experiments and animated video corresponding to three different cases of nonlinear pendulum: a singular pendulum with negligible damping, a heavily damped pendulum, and a chaotic double pendulum. This example differs from the first two in that it focuses on identifying the algebraic constraints and, therefore, a reduced coordinate system, which is not readily apparent from pixel data. For the single-pendulum example, we also demonstrate recovery of a dynamic model in the new coordinate.

In Examples 1 and 2, where the time series were generated synthetically, we asses data requirements, robustness to noise, and sensitivity to different types of perturbations are studied.  For a simulated time-series states, $\mathbf{X}_i$ , Gaussian noise is sampled and added to the  time series as 
\begin{equation}\label{eq:noise}
     \mathbf{X}_{i,\text{noisy}}=\mathbf{X}_i+noise\%\cdot\sigma_i\cdot\mathcal{N}(0,1),  
\end{equation}
where $noise\%$ is the percentage of noise, $\sigma_i=\text{std}(\mathbf{X}_i)$, and $\mathcal{N}(0,1)$ is the standard normal random variable. In Example 3, the pixel data was extracted from video footage, which inherently introduced noise into the dataset. In Examples 1 and 3, we explore the impact of polynomial library size on data requirements.
For both the algebraic finder and the ODE discovery step, we used  LASSO followed by sequential thresholding as the default sparse solver in all our numerical experiments, unless specified otherwise. The values of the hyperparameters representing sparsity enforcement penalties, $\alpha^a$ and $\alpha^d$, for each numerical experiment where applicable, are provided in the Supplementary-section S6.

\subsection{Example 1: Chemical Reaction Networks}\label{Sec:CRN_Example}
We demonstrate the ability of SODAs to accurately recover governing equations from simulated time series of known CRNs, and study data requirements as system complexity and extrinsic noise increase. In applications in chemistry and biology, the temporal evolution of each chemical species in a CRN can be modeled using mass-action kinetics [\citen{erdi1989mathematical,alon2019introduction}]. This naturally leads to models consisting of polynomial terms of the concentrations of chemical species. However, timescale separations in the equilibration rates of different species motivate quasi-steady-state approximations [\citen{bowen1963singular,segel1989quasi}] which simplify mathematical analysis by reducing the state space dimension [\citen{goeke2012computing}].

An advantage of discovering the DAE form rather than the related rational equations is that we discover which variables are in quasi-steady state rather than assuming that all measured variables are inherently dynamic. For known systems, models can be simplified through conservation analysis, partitioned into sub-networks, or non-dimensionalized. Discovering the algebraic systems from data could automate the process of model reduction in chemical and biological networks when the systems are unknown.  Beyond quasi-steady-state approximations, other optimization and sensitivity analysis can motivate simplification [\citen{snowden2017methods}]. Previous work used neural network-based model discovery as a framework for CRN reduction [\citen{SMAI-JCM_2021__7__121_0, benner2020operator}]. SODAs demonstrates identifies likely conservation-of-mass relationships and timescale separation using regression-based sparse-model selection, resulting in an interpretable reduction when model structure is unknown. Note that SODAs requires measurements of both algebraic and differential variables, which is unlikely for realistically-sized biological or chemical networks. Future work to lessen the burden of unobserved states is discussed in Section \ref{Sec:DiscussionSec}. In realistic examples, the method in its current form may be best-suited instead for model reduction during the simulation of large CRNs, which is an ongoing area of study [\citen{yang2016l1,wen2023chemical}].

\subsubsection*{System description}\label{sec:CRN_desc}

\begin{figure}[h!]
    \centering
    \includegraphics[width=1\linewidth]{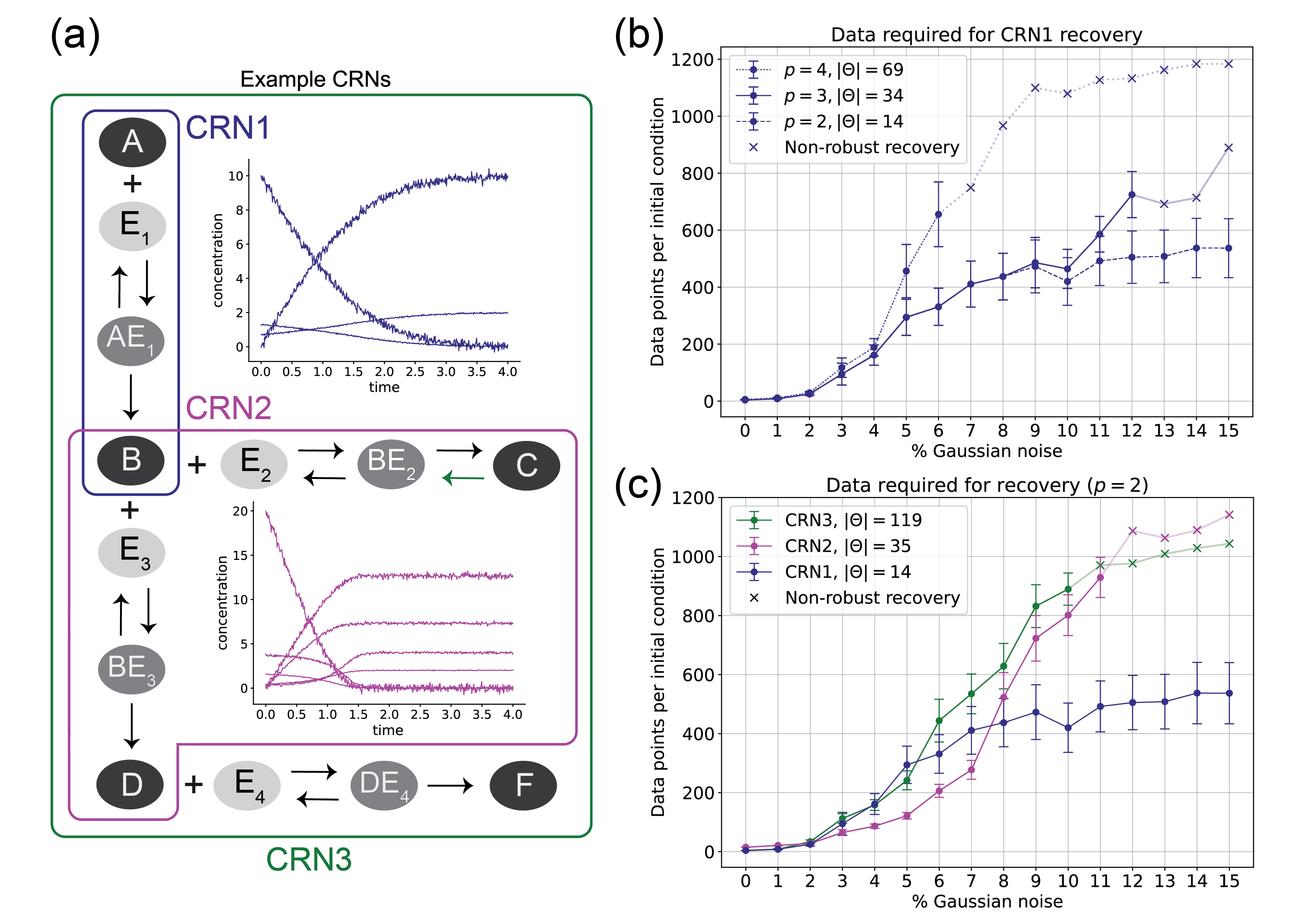}
    \caption{\textbf{Example 1: Application to chemical reaction networks.} (\textbf{a}) Nested structure of CRN1 and CRN2 within CRN3 and examples of simulated time-series with $5\%$ noise. Note that only CRN3 allows reversibility of $C \to BE_2$ as denoted by the green arrow. (\textbf{b}) Scaling of the sampling rate required for recovery of the correct algebraic equations in CRN1 across library degree $p$ and cardinality $|\Theta|$. Non-robust recovery is denoted when fewer than $80\%$ of random seeds succeeded in recovering the correct algebraic relationships. Error bars represent standard error across 10 random noise seeds. (\textbf{c}) Data requirement for recovery of the correct algebraic equations for all CRNs (cardinalities noted) with fixed library degree $p=2$. Each system instance used five initial conditions which are listed in Supplementary-section S3(c).}
\label{fig:CRNs}
\end{figure}

In this example, we use time series generated from three CRNs of varying complexity consisting of one, two, and four enzyme-mediated reactions (Figure~\ref{fig:CRNs}a). The four-reaction system describes the simplified propanediol utilization pathway native to \emph{Salmonella enterica} [\citen{mills2022vertex,kennedy2022linking}]. CRN1 is the single irreversible reaction (\ref{eq:CRN-sim-1}- \ref{eq:CRN-sim-4}.) CRN2 and CRN3 follow a similar mass action kinetic form (full models, parameters, and initial conditions in S3).  CRN1 has two algebraic relationships, one representing conservation of enzyme (\ref{eq:CRN-sim-3}) and one from a quasi-steady state approximation (\ref{eq:CRN-sim-4}). Conservation equations are common in modeling CRNs, representing the notion that enzymes are not created or degraded on short time-scales, and that the system is closed \cite{erdi1989mathematical}. 
CRN2 has four algebraic and three dynamic variables. 
 CRN3 has nine algebraic and five dynamic variables.  Because of the reversibility of the $[B] \Leftrightarrow [C]$ reaction in CRN3, the algebraic relationships have increased complexity. For example, for $[E_2]$, enzyme conservation results in $[E_2]+[B E_2] +[C E_2]-E_2^{tot} =  0$ and the quasi-steady-state approximations $ k_4[B][E_2] - (k_5+k_8)[BE_2] + k_9[CE_2] \approx 0 $.
 Irreversible reactions in CRN2 and CRN3 follow very similar forms to those introduced in CRN1.

 During the simulation of the systems to produce time-series to test SODA's recovery, we assume the algebraic relationships hold exactly using the systems in reduced Michaelis-Menten form (i.e. \ref{eq:CRN-sim-redODE} and \ref{eq:DAE-red-example}, and S3(c)). Each system was simulated from five different initial conditions, varying the initial concentration of the first substrate in each network, so that the range of dynamic behavior was similar in all instances explored.
As SODAs has no knowledge of the difference between algebraic relationships arising from conservation or quasi-steady-state these equations can appear as any linear combination of the equations (See S3(b)).
In the dynamic finder step, we expect to find equations for the dynamics of the substrates. For example in CRN1 this is $A$, and the product $B$, represented by (\ref{eq:CRN-sim-1}) and (\ref{eq:CRN-sim-2}) (See S3(a) for CRN2 and CRN3 systems).

\medskip\noindent
\subsubsection*{SODAs application.} 



We tested SODAs recovery of the expected algebraic relationships from the simulated time series of CRN1, 
using libraries of polynomials with degree $p=2,3,4$ and a range of 0-15\% noise  (Figure \ref{fig:CRNs}b).  The high end of this noise range is similar to realistic noise in metabolite measurements [\citen{bartha2007effect,meyer2012minimum}].  
Data was preprocessed using Savitzky-Golay filters [\citen{savitzky_filter}] to smooth noisy derivatives. Since we were interested in the minimum data to ensure rediscovery for each library and noise condition, we defined the SODAs data requirement to be the smallest number of uniformly spaced data points per initial condition above which the algebraic relationships were recovered for 10 noise instances. The time series generated from five initial conditions were sampled at varied frequencies up to a maximum of 1200 data points per initial condition (see Supplementary-section S3(c) for visualizations of time series). At lower sampling resolution than indicated in Figure \ref{fig:CRNs}b, recovery may be possible but is dependent on the particular realization of measurement error (see Supplementary-section S3(d)). Subsequent discovery of the dynamic equations returned the correct structure and parameter sufficiently similar to the ground truth (see Supplementary-section S3(f) for examples of parameter fitting results at varied noise).

The algebraic relationships were correctly recovered with up to $15\%$ noise for the degrees of $p$ tested for at least some noise instances. For $p= 2$ and up to 11\% noise for $p = 3$, recovery was robust across all noise instances.  The amount of data required for recovery across library degrees scaled similarly for low noise, but the introduction of higher-order terms, especially in the $p=4$ library, led to a sharp increase in the data requirement as noise becomes significant. Due to corruption of the signal and noise amplification, these higher-order terms likely induce spurious algebraic relationships with correlations similar to those of true relationships. For $p=4$, as the noise was increased beyond $6\%$, recovery failed to be robust. There are certain noise realizations for which $\leq1200$ data points are sufficient for rediscovery, but rediscovery is not reliable. SODAs required a lower resolution of data and only five initial conditions for recovery compared to implicit-SINDy [\citen{I_sindy}] or SINDy-PI [\citen{PI_sindy}], which demonstrated recovery using time-series from over 3,000 initial conditions. We also attempted to apply PI-SINDy which is the most robust algorithm among implicit SINDy methods, to discover the models using data where SODAs was successful, but were unable to recover the correct dynamics even when using the maximum amount of data considered in this study (Figure \ref{fig:CRNs}b: 5 initial conditions and 1200 time-series points per initial condition) with as low as $1\%$ noise. The method was applied using the PySINDy implementation of PI-SINDy (details in Supplementary-section S3(g)). However, SODAs requires knowledge of the time-series of complexes and enzymes (ie. $[AE_1]$ and $[E_1]$), while SINDy-PI requires only the timeseries of the substrates (ie. $[A]$), so the comparison is not direct and challenging to interpret.

Testing recovery in the more complex CRN2 and CRN3 examples demonstrates how SODAs responds to a larger state space, richer dynamics, and, as a consequence, more complex algebraic relationships (Figure \ref{fig:CRNs}). Savitzky-Golay filters were again applied [\citen{savitzky_filter}] to smooth noisy derivatives. For CRN2, Pareto front analysis correctly suggested three refinements, while, for CRN3, it correctly suggested five (see Supplementary-section S3(e)). We assessed the method for increasingly noisy data across a variety of resolutions at fixed polynomial degree $p=2$ (Figure \ref{fig:CRNs}c). The data requirement increases with noise, but there is a notable difference in how the data requirement responds to an increase in underlying model complexity vs. an increase in library degree, even if the two modifications lead to $\Theta$ of similar cardinality. For example, $p=2$ recovery for CRN2 requires significantly more data points per initial condition than $p=3$ recovery for CRN1 at $>8\%$ noise, despite the two cases having $|\Theta|=35$ and $|\Theta|=34$, respectively. We hypothesize that SODAs performance is more sensitive to the library degree than the dimension of the library's monomial basis, a notion consistent with our earlier discussions on the impact of higher-order terms on library correlation. Our results also showed that it is not necessarily the case that more complex underlying models always have higher data requirements for rediscovery, as there were regimes where CRN2 required fewer data points per initial condition than CRN1, and regimes where CRN3 required fewer data points than CRN2. This highlights that the quantity of the data alone does not determine success, and the range of dynamics sampled within the data may be important to model discovery.

We finally note that, due to the instability that can arise during forward simulation of discovered DAEs and the inherent difficulty of parameter estimation for DAE systems, we relied solely on the goodness of fit measured by the $R^2$ score as the primary metric for model selection in our numerical experiments. This choice enables a consistent and robust comparison of candidate models without requiring successful forward simulation or stable parameter estimation. In practice, however, a practitioner may instead employ a DAE solver together with a parameter estimation framework tailored to the specific class of DAEs under consideration, using simulation-based validation as an additional criterion for model selection among candidates produced by SODAs under different data regimes, noise realizations, library choices, and hyperparameter settings. We provide a concrete demonstration of this approach in Supplementary-section S3(f), where three representative models discovered at $10\%$  additive noise using a degree-2 library are examined: one with the correct structure and two with incorrect structures of varying similarity to the true model. In this example, parameters are estimated, and the resulting DAEs are simulated and compared against noisy data, illustrating that when the discovered model structure deviates substantially from the true system, the simulated DAEs tend to become unstable, and meaningful parameter inference and simulation become difficult.

\subsection{Example 2: Power-grid Dynamics}\label{Sec:powergridExSec}

The operation of modern power grids relies on detailed modeling of dynamical behavior and network topology, which are critical for grid stability, efficiency, and control against disturbances [\citen{gong2023data}]. Such model information directly supports daily operations, scheduling, and advanced techniques like state estimation, fault monitoring, and optimal power flow computation. However, the rise of renewable generation has led to decentralized power generation, making detailed grid models less accessible in distribution systems In these settings, the power-grid network changes frequently due to the addition or removal of generating units (e.g., rooftop solar, battery storage units) and network components (e.g., lines, switches, transformers), as well as uncertainty in distributed energy resources, line maintenance, faults, and other unplanned events. Yet, traditional power-grid monitoring and fault detection methods are often grounded on parameter estimation techniques, requiring extensive prior knowledge of grid configurations and physical properties [\citen{pierre2012overview,gopakumar2018remote,liu2021dynamic,gaskin2024inferring}]. The increasing deployment of phasor measurement units (PMUs) [\citen{rodrigues2021pmu}], capable of providing high-resolution data, offers a significant opportunity for automated model discovery [\citen{ardakanian2019identification,yuan2019data}]. Here, we show that the sparse discovery of differential-algebraic models can robustly infer the network topology, which is encoded by the algebraic relationships among power injections at each bus. Using numerical results, we demonstrate that the proposed method accurately infers network models across a wide range of scenarios, including both small- and large-scale power grids, varying levels of state perturbation, and diverse PMU signal qualities.

\medskip\noindent
\subsubsection*{System description.} The dynamics of power grids can be modeled as a network of coupled oscillators (nodes) interconnected through power transmission lines (links) [\citen{sadamoto2019dynamic,dorfler2013synchronization,nishikawa2015comparative}]. The nodes correspond to conventional generators and consumers (geographically represented by substations and their attached load). In this representation, the dynamics of generators are governed by the {swing equation} [\citen{saadat1999power}]: 
\begin{equation}
\frac{2H_i}{\omega_{\rm R}}\ddot{\phi}_i + \frac{D_i}{\omega_{\rm R}}\dot{\phi}_i = P_{\rm m, i} - P_{\rm e, i}, \,\,\, \text{for} \,\, i=1,\ldots,n_{\rm g},
\label{eq.generator}
\end{equation}
and the dynamics at the load buses and generator terminals follow the first-order phase oscillator model: 
\begin{equation}
\frac{D_i}{\omega_{\rm R}}\dot{\phi}_i = P_{\ell,i} - P_{\rm e, i}, \,\,\, \text{for} \,\, i=n_{\rm g}+1,\ldots,N,
\label{eq.load}
\end{equation}
where $n_{\rm g}$ is the number of generators, $n_{\ell}$ is the number of load buses, and $N = 2n_{\rm g}+n_{\ell}$ is the number of oscillators. The generator terminals correspond to buses at which generators are linked to the transmission system, acting as boundary points where the electrical output of the generators is delivered to the grid. 
Here, $\phi_i(t)$ is the phase angle of oscillator $i$ at time $t$ relative to the frame rotating at reference frequency $\omega_{\rm R} = 2\pi f_{\rm R}$, while $H_i$ and $D_i$ are the inertia and damping constants accounting for the mechanical properties of the generator. The supplied mechanical power by generators is $P_{{\rm m},i}>0$, and the consumed power by loads is $P_{\ell,i}<0$. For each bus $i=1,\ldots,N$, the electrical active power $P_{{\rm e},i}$ and reactive power $Q_{{\rm e},i}$ are determined by the system's conservation law (power flow) according to the following algebraic equations:
\begin{equation}
\begin{aligned}
        P_{{\rm e},i} &= +\sum_{j=1}^N | V_i V_j Y_{ij} | \sin (\phi_i - \phi_j - \gamma_{ij}), \\
        Q_{{\rm e},i} &= - \sum_{j=1}^N | V_i V_j Y_{ij} | \cos (\phi_i - \phi_j - \gamma_{ij}),
\end{aligned}
\label{eq.algeqs}
\end{equation}north

\noindent
where $V_i$ is the line voltage at bus $i$ (assumed to be constant). The matrix entry $Y_{ij} = |Y_{ij}|e^{i(\gamma_{ij}+\frac{\pi}{2})}$ represents the admittance of the transmission line interconnecting buses $(i,j)$ in polar form. Note that, at steady state, the system is balanced (i.e., $\sum_{i}P_{\rm e,i} = 0$). For simplicity, we assume that the line's impedances are purely inductive (i.e., $\gamma_{ij} = 0$).

\medskip\noindent
\subsubsection*{SODAs application.}
\begin{figure}[!ht]
    \centering
    \includegraphics[width=1\linewidth]{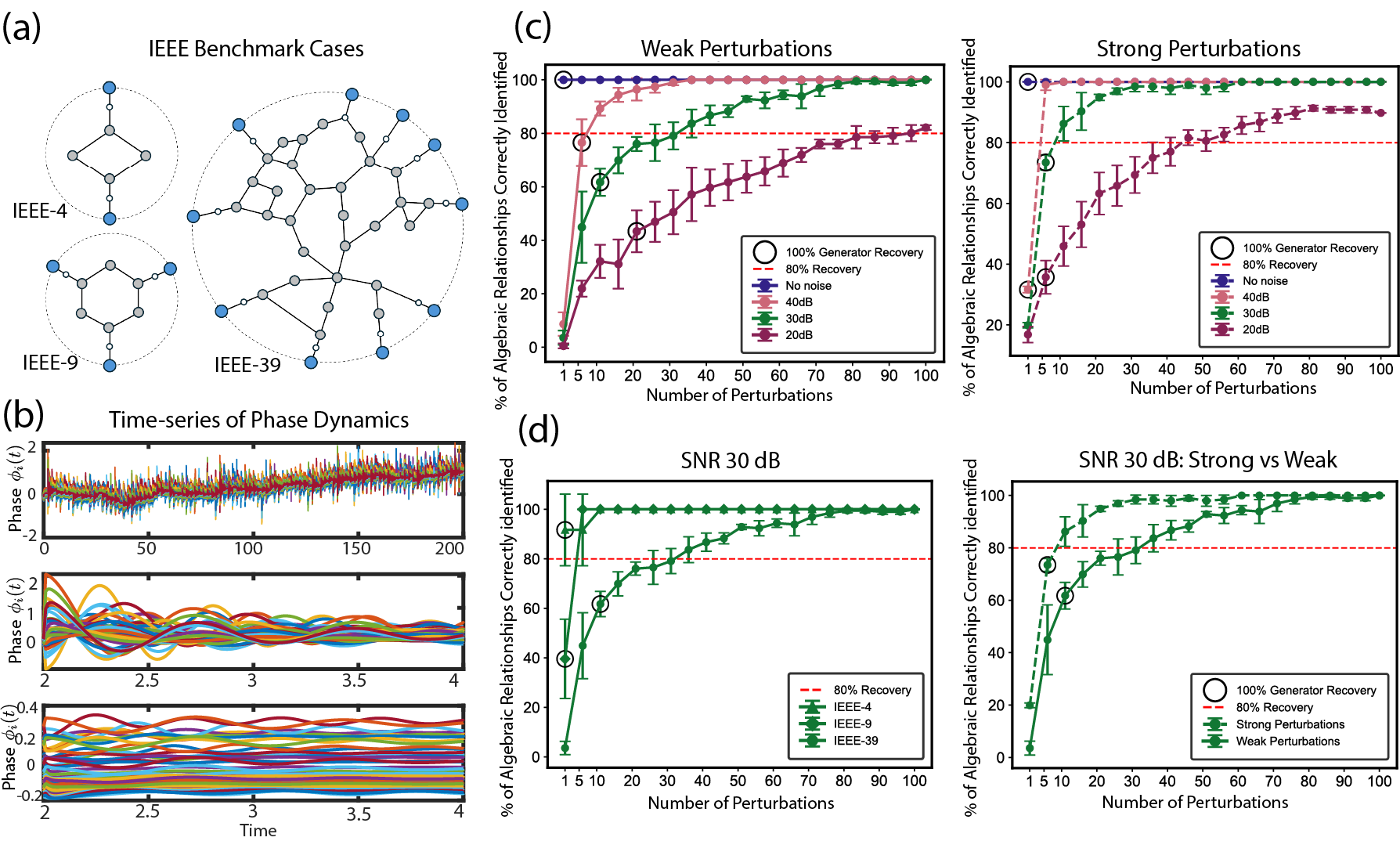}
    \caption{\textbf{Example 2: Discovery of power-grid networks.}
    (\textbf{a}) Power-grid networks of the IEEE-4 (\textit{top left}), IEEE-9 (\textit{bottom left}), and IEEE-39 (\textit{right)} benchmark systems. Generator buses are represented by blue nodes on the outer circle, while aggregated load buses are represented by gray nodes. 
    (\textbf{b}) Time series of phase dynamics $\phi_i(t)$ for all nodes in the IEEE-39 power grid under different levels of perturbations: large perturbations (\textit{top}), large perturbations zoomed-in (\textit{middle}), and small perturbations zoomed-in (\textit{bottom}).
    (\textbf{c}) Performance comparison of SODAs  applied to IEEE-39 benchmark system across different SNR levels and : small perturbations (\textit{left}) and large perturbations (\textit{right}).
    (\textbf{d}) Performance comparison of SODAs  at 30 dB SNR across: different benchmark cases (\textit{left}) and different perturbation strengths (\textit{right)}.}
    
    \label{fig:Example_3_powergrid_fig}
\end{figure}

We considered three IEEE benchmark cases with varying complexity (Figure \ref{fig:Example_3_powergrid_fig}a): IEEE-4, IEEE-9, and IEEE-39. Time series were generated under two distinct conditions: i) strong perturbations (Figure \ref{fig:Example_3_powergrid_fig}b - \textit{top}), which excite nonlinear transient behaviors, and ii) weak perturbations (Figure \ref{fig:Example_3_powergrid_fig}b - \textit{bottom}), suitable for small-signal stability analysis. We assume that the power grid is equipped with phasor measurement units (PMUs) placed on all $N$ buses, providing us with high-resolution, high-frequency data for phase angle $\phi_i(t)$, frequency $\dot\phi_i(t)$, line voltage $V_i(t)$, and electrical power $P_{\rm e,i}(t)$ and $Q_{\rm e,i}(t)$. More information on network details, simulations, measurement frequency,   and practicality of assumptions can be found in the Supplementary-section S4(a). For the algebraic finder step, STOLS was utilized instead of LASSO followed by thresholding, as the latter indicated a poor fit when initially applied.



SODAs was applied to the simulated time series to identify the network structure and parameters of the governing DAE system as described in equations \eqref{eq.generator}--\eqref{eq.algeqs}. By using the well-established physical principles of power-grid dynamics, we achieved a significant advantage in constraining SODAs within a more informed framework. More specifically, the conservation of active power at each node leads to inherent algebraic constraints that we exploited to improve the identification process. In constructing the candidate library, we included the following functional forms: $P_{{\rm e},i}$ (active power states), $\phi_i$ (phase angles), $\dot\phi_i$(frequency), and sine terms of phase differences $\sin(\phi_i - \phi_j)$ for all pairs $(i,j)$. In a network with $N$ nodes, this initial library included $3N + \frac{N(N+1)}{2}$ terms, which amounts to 39, 114, and 1,372 terms for the IEEE-4, IEEE-9, and IEEE-39 benchmarks, respectively. Additionally, by incorporating prior knowledge, we further reduced the number of candidate terms used to discover each conservation law. For example, the conservation law involving $ P_{{\rm e},i}$ at a given node $i$ restricts the algebraic relations to include only sine terms of the form $\sin(\phi_i - \phi_j)$, for all nodes $j$. This restriction decreases the number of coupling terms for each node from $\frac{N(N+1)}{2}$ to $N-1$. Consequently, the number of library terms used to discover each power-flow equation consisted of 23, 47, and 195 terms for the IEEE-4, IEEE-9, and IEEE-39 benchmark cases, respectively.
Using the \texttt{Daefiender.FeatureCouplingTransformer} class, we efficiently enforced these additional library constraints, minimizing computational complexity and enhancing the accuracy of identified models by eliminating irrelevant terms.

In the algebraic finder step, conservation laws constraining the active powers $P_{{\rm e},i}$ for all $i$ are discovered \eqref{eq.algeqs}. These algebraic equations inherently reveal the power grid's network structure (encoded by the admittance matrix $Y$). In the dynamic finder step of SODAs, the differential equations governing the generators are discovered using a SINDy framework. We approximated the left-hand side of the dynamic equations \eqref{eq.generator}--\eqref{eq.load} by utilizing derivatives of spline-fitted phase angles and applied sparse regression techniques, where LASSO was followed by sequential thresholding. Parameter estimation was subsequently performed through linear fitting with the discovered model structure. Note that for the identification of the dynamic equation, we modified the candidate library by removing all \( \sin(\phi_i - \phi_j) \) terms, which had already been accounted for by the algebraic constraints. Instead, we incorporated state measurements of \( P_{{\rm e},i} \) and hence reducing the number of candidate functions in the dynamic finder step to $3N$. This refinement has also helped to avoid using the coupling terms \( \sin(\phi_i - \phi_j) \), which have high sensitivity to noise in the data, for the dynamic finder step. We observed complete recovery of the true ODE system whenever the algebraic constraints were discovered with 100\% accuracy. In contrast, attempts to identify coupled ODE equations without segregating the algebraic terms led to poor model recovery, even at SNR as high as 40dB, emphasizing the importance of this methodological separation of algebraic and differential equations. 

We present SODAs' performance on identifying the ground-truth model for all three test cases, with a primary focus on the IEEE-39 benchmark test case (Figure \ref{fig:Example_3_powergrid_fig}c and \ref{fig:Example_3_powergrid_fig}d). We conducted similar analyses on the IEEE-4 and IEEE-9 benchmark networks, where SODAs demonstrated performance comparable or superior to that of the more complex IEEE-39 case. To evaluate the robustness of SODAs, we introduced varying amounts of noise into the simulated time-series data and examined performance across four different signal-to-noise ratios (SNR): no noise, 40 dB, 30 dB, and 20 dB. For each SNR level, we generated five distinct noise realizations and analyzed the performance as an average across these realizations. We found no significant variations in performance among the different noise realizations within the same SNR level, which allowed us to exclude error bars related to noise realizations in our figures.

The dynamics of the system are represented in the time series by the transients created from perturbations. Consequently, the number of perturbations reflects the amount of transient information about the system's dynamics. To evaluate SODAs' performance on this time-series dataset, we measured the percentage of true algebraic relationships recovered as a function of the number of perturbations in the time series. Additionally, to account for the variability in how the power grid responds to different random perturbations, we used ten different permutations of the perturbation order in our time-series simulations. We averaged the recovery performance over these ten permutations, representing the associated variability with error bars in our results. 

For SNRs of no noise, 40 dB, and 30 dB, SODAs achieved 100\% recovery of all algebraic relationships, as well as the dynamic ODEs, using the refined library for both small and large perturbations. Even under high noise conditions (SNR of 20 dB), SODAs maintained a recovery rate exceeding 80\% for algebraic relations with 100 perturbations. The algebraic relationships related to power conservation at generator nodes were identified with relatively small number of perturbations, underscoring the significant role of these nodes in influencing the overall dynamics of the power grid. Techniques such as signal averaging and advanced filtering can improve the SNR, even when the inherent measurement SNR is lower. Therefore, an SNR of 30 dB is assumed to be realistic for our further analysis, as other studies have reported PMU devices with signal quality at or above this level [\citen{brown2016characterizing,bhandari2019real,frigo2019statistical,zhang2020experimental}]. We evaluated the performance of SODAs at a 30 dB SNR across various benchmark cases and levels of perturbations (Figure \ref{fig:Example_3_powergrid_fig}d). We did a comparative analysis of the three benchmark cases, in terms of the number of perturbations required to recover all true DAE system components, and found that the number of perturbations required for full recovery increased with the complexity of the power grid, yet it remained within the realistic number of permutations observed in practice (Figure \ref{fig:Example_3_powergrid_fig}d - \textit{left}). Also, strong perturbations that induce significant transients in system dynamics enabled quicker recovery of the DAE system as expected (Figure \ref{fig:Example_3_powergrid_fig}d - \textit{right}). 

Additionally, we also investigated whether the models discovered by SODAs produce physically meaningful power system trajectories when simulated. Across multiple test cases, we found that accurately recovering the algebraic (power-flow) constraints was crucial for generating valid phase angle dynamics. When these algebraic relationships were misidentified, the resulting DAE models often failed to demonstrate frequency synchronization. This failure was evident as the phase angles did not converge to a steady equilibrium in the co-rotating frame, indicating that the predicted system dynamics were qualitatively incorrect. A representative example from the IEEE-4 benchmark system, illustrating the differences between correctly and incorrectly identified algebraic constraints, can be found in Supplementary-section S4(b).

\subsection{Example 3: Non-linear and Chaotic Pendulums}\label{Sec:pend_ex}
In this example, we demonstrate how the algebraic finder step of SODAs can identify a reduced coordinate system, specifically the polar coordinate system that arises from the physical constraints of three different pendulums. This example differs from the previous two examples, where the goal was to derive a complete set of DAEs that govern the system’s dynamics. In this example, we apply SODAs to analyze pixel data obtained from video recordings of three different types of pendulums:  1) real-time footage of a single pendulum exhibiting non-linear behavior through large-angle swings, 2) video recordings of the animation of a highly damped pendulum, and 3) video recordings of the animation of a double pendulum exhibiting chaotic behavior. Screen recordings of the animation were used for cases 2 and 3 to mimic the process of extracting data from video, even though a simulation was used to produce the animation. These cases illustrate the use of SODAs in mechanical systems, where state measurements are obtained through images, videos, or sensor data.

When modeling the motion of mechanical systems, there often exist coordinate transformations that lead to simple equations explaining the motion. For example, the natural coordinate system for modeling the motion of a single pendulum is known to be the one-dimensional polar coordinate system, i.e., to model the angle $\omega$ that the pendulum makes with a reference axis in the 2D plane (Figure \ref{fig:Example_2_pendulum_fig}a - \textit{top}) or in the case of a double pendulum, the two angles, $\omega_1,\omega_2$ , that the two pendulums make with the reference axis (Figure \ref{fig:Example_2_pendulum_fig}a - \textit{bottom}). However, finding the appropriate coordinate transformation to produce the most parsimonious system is often challenging and not clear from raw measurements of the states. Recently, several data-driven methods have been used to find reduced coordinates, for example, using neural networks to learn nonlinear mappings to a lower-dimensional state-space [\citen{sindy_auto_2024bayesian,yang2023latent}]. Here, we show how SODAs can assist in this process using the algebraic finder step. Specifically, discovering the algebraic relationships between state variables can give intuition about the nature of symmetries within the system and the natural coordinate system suitable for obtaining a simple model. 

In all three cases presented here, SODAs identify the intrinsic algebraic constraints in the system, implying the reduced polar coordinate system inherent to it, solely based on pixel data that represents the time series of the pendulum's position in Cartesian coordinates $(x, y)$. For Cases 1 and 2, once the coordinate transformation is identified, the ODE equation of the pendulum is discovered in the new transformed coordinate system. For Case 3, we do not recover the dynamic equations. Together, these cases highlight that discovering algebraic constraints can aid in finding a coordinate system and is less data-intensive than discovering the complete DAE system.

\medskip\noindent
\subsubsection*{System description.}
A non-linear single pendulum with a mass attached to its end can be modeled using the following second-order ODE:

\begin{equation} \label{eq:nonl_pend}
m\ddot{\omega}+\alpha\dot{\omega} +m \frac{g}{l} \sin \omega  = 0,
\end{equation}
where $m$ and $l$ represent the mass and length of the pendulum, respectively, $\omega$ denotes the angle that the pendulum makes with the vertical axis, $g$ represents the acceleration due to gravity, and $\alpha$ is the damping coefficient, which is zero when damping is negligible. The physical constraint of the single pendulum can be captured by the equation
\begin{align}
&x^2 + y^2 = l^2, \label{eq:pend_constraint}
\end{align}
where the transformations between the polar and Cartesian coordinates are given by:
\begin{align}
&x=l\sin \omega;\quad  y=-l\cos \omega .\label{eq:polar_transf}
\end{align}

\begin{figure}[!ht]
    \centering
    \includegraphics[width=1\linewidth]{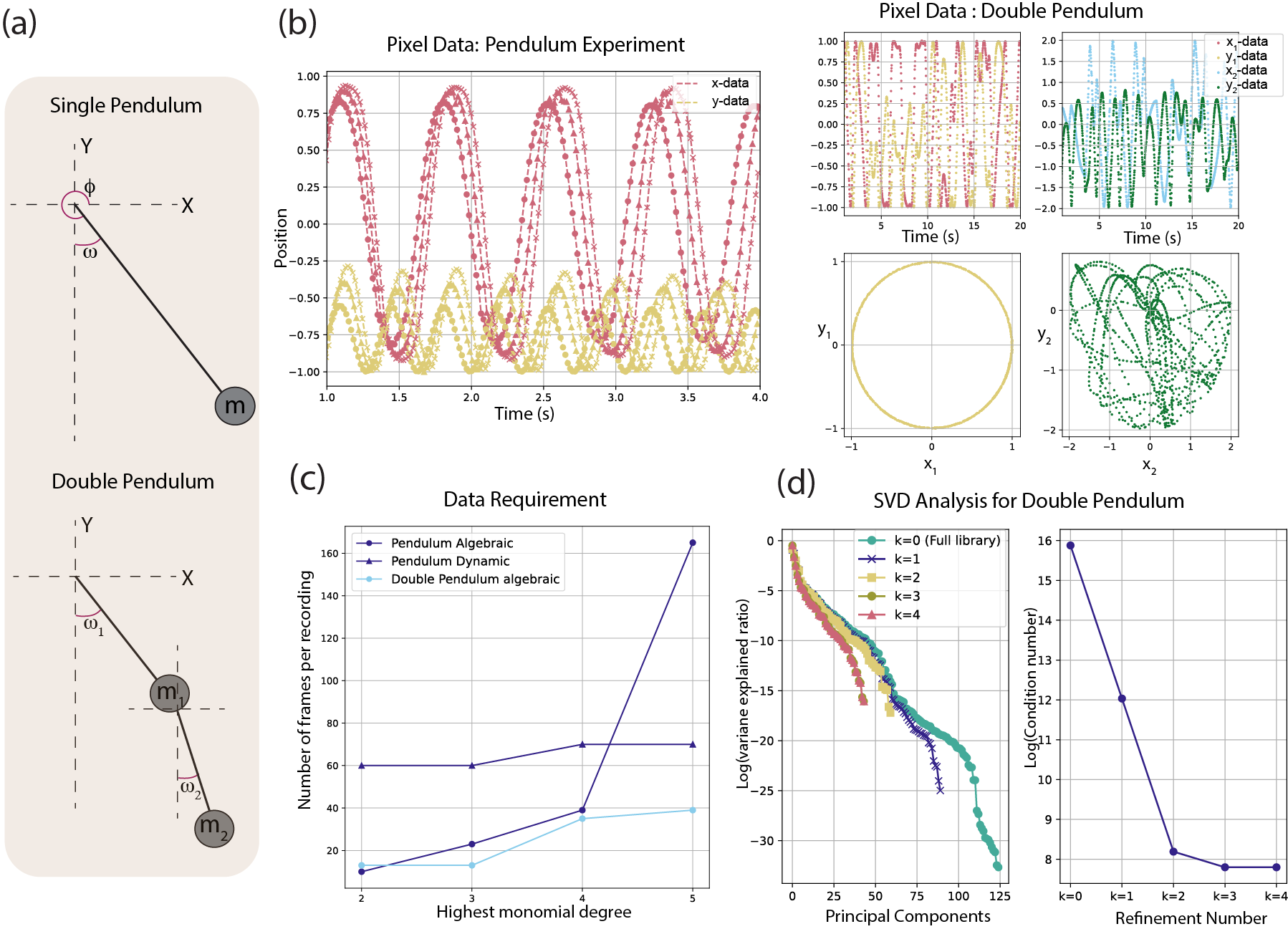}
    \caption{\textbf{Example 3: Application to non-linear pendulum.}
       (\textbf{a}) Schematic diagram: single pendulum (\textit{top}) and double pendulum (\textit{bottom}).  
       (\textbf{b})  Scaled pixel data: single pendulum experiments (\textit{left}) and  double pendulum animation (\textit{right}). 
    (\textbf{c}) Data requirement: comparison between damped single pendulum and double pendulum.
    (\textbf{d})  SVD analysis: determining the number of algebraic constraints in the case of the  chaotic double pendulum with degree 5 library. Y-axes are given in natural log scale.}
    
    \label{fig:Example_2_pendulum_fig}
\end{figure}

A double pendulum consists of two pendulums attached end to end, as shown in Figure \ref{fig:Example_2_pendulum_fig}a - \textit{botom}. This assembly of two pendulums is well-studied and known to exhibit chaotic motion [\citen{shinbrot1992chaos}]. To simplify the analysis, we will neglect damping. The equations of motion for a double pendulum, derived using Lagrangian mechanics [\citen{shinbrot1992chaos}], are given by: 
\begin{align}
&\ddot{\omega}_1 = \frac{g (\sin \omega_2 \cos (\omega_1-\omega_2) - \mu \sin \omega_1) - (l_2 \dot{\omega}_2^2 + l_1 \dot{\omega}_1^2 \cos (\omega_1-\omega_2)) \sin (\omega_1-\omega_2)}{l_1 (\mu - \cos^2 ( \omega_1-\omega_2))}, \label{eq:double_pend_eq_1} \\
&\ddot{\omega}_2 = \frac{g \mu (\sin \omega_1 \cos (\omega_1-\omega_2) - \sin \omega_2) + (\mu l_1 \dot{\omega}_1^2 + l_2 \dot{\omega}_2^2 \cos (\omega_1-\omega_2)) \sin (\omega_1-\omega_2)}{l_2 (\mu - \cos^2 (\omega_1-\omega_2))}, \label{eq:double_pend_eq_2} 
\end{align}
where $\omega_1$ and $\omega_2$ denote the angles that the first and second pendulums make with the vertical axis, respectively, $l_1,l_2$ denote their respective lengths, $m_1,m_2$ denote their respective masses, and $\mu=1+\frac{m_1}{m_2}$. The physical constraints in this mechanical system are given by 
\begin{align}
&x_1^2 + y_1^2 = l_1^2, \label{eq:double_pendulum_constr_1}\\
&(x_1-x_2)^2 + (y_1-y_2)^2 = l_2^2,\label{eq:double_pendulum_constr_2} 
\end{align}
where $(x_1,y_1)$ and $(x_2,y_2)$ denote the positions of the first and second pendulum in Cartesian coordinates. The transformations between the Cartesian and polar coordinate systems in this case are given by $x_1=l_1\sin \omega$, $y_1=-l_1\cos \omega_1$,  $x_2-x_1=l_2\sin \omega_2$, and $y_2-y_1=-l_2\cos \omega_2$.  

\subsubsection*{SODAs application.}
For this study, pixel data from three different cases were used: 1) \textit{single pendulum with negligible damping}: video recording of single pendulum experiments with negligible damping and three different initial conditions (Figure \ref{fig:Example_2_pendulum_fig}b - \textit{left}), 2) \textit{heavily damped pendulum}: video recording of 
 animation of a heavily damped pendulum with four different initial conditions, 3) \textit{chaotic double pendulum}: single video recording of the animation of a double pendulum exhibiting chaotic behavior (Figure \ref{fig:Example_2_pendulum_fig}b - \textit{right}). OpenCV [\citen{openc_cv_book}] package in Python was used to extract the pixel data from the video recordings. We would like to note that converting raw video or animation frames into reliable image-plane measurements is a complex problem. This process is managed using a standard, well-established computer vision routine implemented with OpenCV’s standard functions for color segmentation, morphology, and contour detection, [\citen{openc_cv_book,pulli2012real}].  Managing factors such as camera perspective, lens distortion, focus, and noise associated with the video footage is a nontrivial challenge, and hence we relied on these established methods during this stage. More details about the collection of the pixel data and preprocessing can be found in Supplementary-section S5(a).

An initial candidate library was constructed using monomials in the original states: $\Theta_1=\{1,x,y,x^2,xy,\ldots\}$ for Cases 1 and 2, and $ \Theta_1=\{1,x_1,y_1,x_2,y_2,x_1^2,\ldots\}$ for Case 3. These choices were made under the assumption of lack of prior knowledge about the algebraic constraints present in the pixel data and further assuming that these constraints could be captured by a polynomial functions of the state variables. For the case of real-time experimental footage for case 1, monomials up to degree $3$ were included in the candidate library. For case 2 and case 3 where we have better sampling of the state space compared to case 1, the highest monomial degree was varied between $2$ and $5$. A comparative study on the data requirement for equation discovery was carried out for the latter two cases (Figure \ref{fig:Example_2_pendulum_fig}c).

SODAs successfully discovered the constraint \eqref{eq:pend_constraint} for the single pendulum Cases 1 and 2. Note that $l=1$ for our dataset due to the scaling of the pixel data. After discovering the constraint and identifying the reduced polar coordinate system, we applied the following transformation to move to the polar coordinate system: 
\begin{align}
&x=l\cos \phi;\quad  y=l\sin\phi , \label{eq:pendulum_new_polar_transf}
\end{align}
where $\phi$ is the angle formed by the mass of the pendulum with the positive $x$ -axis (Figure \ref{fig:Example_2_pendulum_fig}a - \textit{top}). We made this choice instead of the natural one in\eqref{eq:polar_transf} because, given \eqref{eq:pend_constraint}, transformation \eqref{eq:pendulum_new_polar_transf} is the most straightforward coordinate transformation in the absence of additional information about the system. Indeed, under the translation $\phi=\omega + \frac{3\pi}{2}$, it can be easily verified that the equations we discovered after applying the transformation \eqref{eq:pendulum_new_polar_transf} are consistent with \eqref{eq:polar_transf} and \eqref{eq:nonl_pend}.
 For this simple case, we have sufficient data to discover the dynamic equations in the polar coordinate system. We applied a Savitzky-Golay filter [\citen{savitzky_filter}] with a window size of 12 to compute the first and second derivatives of the $\phi$. The candidate library consisted of monomials of $\phi,\dot{\phi},\cos \phi$ and $\sin \phi$ : $\Theta_*=\{\phi,\dot{\phi}, \cos \phi,\sin \phi,\phi^2,\phi \dot{\phi},\ldots \}$. Using \eqref{eq:pendulum_new_polar_transf}, we included $x$ and $y$ in the library instead of calculating $\cos \phi$ and $\sin \phi $ directly from $\phi$. Similar to the algebraic finder step, we restricted the monomial degree to 3 for Case 1 and varied it between $2$ and $5$ for Case 2 for a comparative study.
 After ODE model structures were discovered using sequential LASSO and thresholding, we used a curve fit package to determine the correct parameters. In both case 1 and case 2, simulations of the discovered equations predicted the test pixel data accurately (see Supplementary-section S5(b) for details on discovered equations and parameter estimation). Noteably, the discovery of the dynamic equations in the library with $\sin\phi$ and $\cos\phi$ required more data than the algebraic equations in the original Cartesian coordinates for monomial degree 5 (Fig. \ref{fig:Example_2_pendulum_fig}c), likely due to a combination of higher correlations in the candidate library and the need to calculate derivatives from noisy data for ODE discovery.

For Case 3, SODAs robustly identified the constraints \eqref{eq:double_pendulum_constr_1}--\eqref{eq:double_pendulum_constr_2} from pixel data obtained from a single video segment, using various upper bounds on the degree of the monomial library. Note that the polar nature of the states $(x_2,y_2)$ is not immediately evident from the pixel data of the secondary pendulum in the double-pendulum (Figure \ref{fig:Example_2_pendulum_fig}b - \textit{right}). Nevertheless, SODAs discovered an expanded version of \eqref{eq:double_pendulum_constr_2}. This algebraic discovery enables us to determine the transformation to the polar coordinate system: $\phi_1=\arctan \frac{y_1}{x_1}$ and $\phi_2=\arctan \frac{y_1-y_2}{x_1-x_2}$, which could, in principle, be further used for dynamic discovery. The dynamic discovery of equations governing the motion of the double pendulum similar to the form \eqref{eq:double_pend_eq_1}
--\eqref{eq:double_pend_eq_2} based on polar coordinates, has been successfully demonstrated using SINDy-based methods [\citen{PI_sindy}]. However, discovering the algebraic constraint from a particular data set does not guarantee the successful recovery of the ODE system in the new coordinate system, and previously successful SINDy-based methods required more data than used here. Success depends on several factors, including the information content in the data relative to the dynamics of the system, noise-robust estimation of derivatives, and more. We suspect discovering this ODE requires weak approaches to handle the higher-order derivatives of \(\phi_1\) and \(\phi_2\). Without such methods, we are unable to rediscover the ODE system defined by equations \eqref{eq:double_pend_eq_1}–\eqref{eq:double_pend_eq_2} as part of this example, but note that even when data is insufficient to discover the dynamics, it can still be used to learn about constraints. To demonstrate the robustness of algebraic equation discovery, we show that the SVD analysis (Section \ref{Sec:Methods}) can be used to evaluate how many library refinements lead to meaningful equations.

We used a Pareto front of the natural log of the condition number of the library to determine the stopping criterion. In Case 3, we stopped after two iterations as no further improvement in the condition number was observed from the third refinement onward (Figure \ref{fig:Example_2_pendulum_fig}d), implying the existence of two algebraic constraints. We also conducted a comparative study of the data requirement for Cases 2 and 3 as a function of maximum degree of the monomials in the candidate library (Figure \ref{fig:Example_2_pendulum_fig}c). The results indicate that the data requirement, in terms of the number of frames per video segment, stayed within manageable limits. Interestingly, the double pendulum required fewer frames per video segment for successful recovery of the algebraic equations compared to Case 2. This is because the chaotic nature of the pendulum, combined with a 3-minute-long video, provided a good sampling of the state space. We also shuffled the data points along the time axis only for Case 3 to ensure a good representation of the dynamics, given the longer length of the video segment. These factors contributed to the reduced data requirement in terms of the number of frames needed for successful recovery in Case 3. An alternative approach to reduce the number of frames required would be to conduct experiments with more informative initial conditions, which is closely related to experimental design [\citen{seltman2012experimental,versyck1999introducing}].

\section{Discussion}\label{Sec:DiscussionSec}

We have introduced SODAs as a method to discover DAEs from data in their explicit form, where the algebraic and differential components are identified separately. This separation played a crucial role in addressing the challenge of perfect multicollinearity within the candidate library, which often arises due to inbuilt algebraic relationships in a DAE system. Through an iterative refinement process, SODAs effectively reduces multicollinearity, ensuring robust and accurate discovery of governing equations. Separation of the algebraic and differential components of the DAE system also ensured that the potential noisy derivatives in the dataset do not affect the discovery of the algebraic component. Crucially, we have demonstrated the flexibility of SODAs by applying the method to three distinct applications. Our results show that SODAs is robust to noise and perform better than other SINDy-based implicit methods. We have also demonstrated that the data requirements for SODAs scale efficiently with the number of data points and the complexity of the library terms used.

There are several limitations and challenges worth pursuing in future work. One significant limitation of SODAs is the necessity to observe or measure all the state variables. In many physical systems, measuring all state variables can be challenging due to cost or physical infeasibility. Previous methods such as implicit-SINDy [\citen{I_sindy}] and SINDy-PI [\citen{PI_sindy}] required no measurements of the algebraic states, but also assumed that all states were differential variables. Extending SODAs to discover DAEs even in the presence of unobserved states, potentially through recent strategies [\citen{DAHSI_paper,sindy_latent,stepaniants2024discovering}], would enable the identification of algebraic and differential variables and their equations when only a subset of states are measured. However, recovery in this setting will be even more challenging due to the non-convexity of sparse-optimization with hidden states as well as the potential lack of analytical identifiability in nonlinear DAE models [\citen{montanari2024identifiability}].

One aspect of DAE discovery we did not yet explore is how SODAs performs when "fast" states are not fully in quasi-steady state. In realistic experimental sampling of systems with fast-slow variables, some transient dynamics of the fast variables may be sampled before they reach quasi-steady state. In this work, we used the reduced ODE formulation \eqref{eq:DAE-red-example} of our CRNs to simulate the time series used for equation re-discovery, ensuring that the slow-variables were in quasi-steady state and thus, there was no error in the algebraic relationships. Future work is needed to investigate how SODAs performs as a function of sampling transient vs quasi-steady-state regimes for fast variables. An interesting avenue would be in mimicking the hybrid-SINDY approach [\citen{mangan_hybrid_sindy,li2023discover}] in which a different model is selected when variables are in quasi-steady state.

Additionally, the discovery of conservation laws for the CRNs could be explored further. Because the total substrate, $M^{tot}$, varied between simulations, identification of the substrate conservation relationship was not possible in our setup, as a single learned $M^{tot}$ cannot fit all simulations. The enzyme conservation relationship could be recovered, as $E_\text{tot}$ was held constant between simulations. Inclusion of non-zero initial conditions in $\Theta$ would allow for the discovery of substrate conservation, a close analog to the inclusion of suspected bifurcation parameters in $\Theta$ in the original formulation of SINDy [\citen{sindy}] and following work. We chose here to exclude initial conditions from $\Theta$ because when any two (or more) initial conditions are held constant across simulations, they will be in perfect correlation, confounding one another and limiting the interpretability of any resulting algebraic relationships. Additionally, discovery of conservation relationships of this nature up to a constant (e.g., relationships mirroring $\frac{d}{dt}([AE_1] + [E_1]) = 0$) was considered, it was not included algorithmically. Inclusion of state derivatives in $\Theta$ would allow for the discovery of relationships of this type, but would require the smoothing of noisy derivatives. We suspect integration with a weak approach \cite{W_sindy} may be more appropriate for conservation laws up to a constant.

There is a more general data limitation challenge inherent to sampling dynamical systems: sampling is often determined by the trajectories of the dynamics themselves which induces correlation between library terms. In cases where the available data does not adequately differentiate between lower-order and higher-order terms, strong correlations separate from the true algebraic relationships can occur. In Example 3, we showed how the data requirement increases as the degree of the monomial terms in the library increases to 5. This is attributed to the high correlation between higher-degree monomials in the library. This high correlation between library terms makes it challenging to separately identify their contributions to the dynamics of the system. This issue is not exclusive to SODAs; it is an inherent challenge in all SINDy-type methods that rely on sparse regression with function libraries. This challenge can be overcome by implementing improved data collection strategies or experimental design techniques [\citen{seltman2012experimental,versyck1999introducing}] that help with distinguishing different candidate library terms. We aim to develop methods to address these challenges in our future work. Finally, this study has focused on DAEs involving only ordinary derivatives, an exciting avenue for future research will be to extend SODAs to DAEs involving partial derivatives. 

\medskip\noindent
\textbf{Data Accessibility: } 
All code and data used for numerical experiments in this paper are available at
\newline
\href{https://github.com/mjayadharan/DAE-FINDER_dev/tree/main/Examples}{https://github.com/mjayadharan/DAE-FINDER\_dev/tree/main/Examples}. 
\newline
The \texttt{DaeFinder} package is also available through Python Package Index (PyPI) at \href{https://pypi.org/project/DaeFinder}{\\https://pypi.org/project/DaeFinder}. 

\medskip\noindent
\textbf{Authors' contributions: } 
M.J. developed the SODAs algorithm, implemented the \texttt{DaeFinder} package, and conducted all experiments (except for Example 1). C.C. contributed to the \texttt{DaeFinder} package and designed and conducted the numerical experiments in Example 1. A.N.M. designed the test cases and assisted in the numerical experiments in Example 2. N.M.M. designed and supervised the research. All authors edited the manuscript and approved the final draft for publication.


\medskip\noindent
\textbf{Competing Interests: } 
We declare we have no competing interests. 

\medskip\noindent
\textbf{Funding: } 
This material is based on work supported by several funding agencies. N. M. M. and M. J. were supported by the U.S. Department of Energy, Office of Science, Office of Advanced Scientific Computing Research, under Award Number DE-SC0024253 and Army Contracting Command under Award Number W52P1J-21-9-3023. C. C. was supported by the National Science Foundation under Grant No. DGS-2021900. Any opinions, findings, and conclusions or recommendations expressed in this material are those of the author(s) and do not necessarily reflect the views of the Department of Energy, Army Contracting Command, or National Science Foundation.

\medskip\noindent
\textbf{Acknowledgments: } 
We thank Finn Hagerty for his assistance in running simulations for the CRN example and related discussions.

\bibliography{sodas}

\end{document}